\documentclass[a4paper]{article}
\pdfoutput=1

\usepackage[english]{babel}
\usepackage[utf8]{inputenc}
\usepackage[colorinlistoftodos]{todonotes}
\usepackage{amsmath}
\usepackage{amssymb}
\usepackage{amsthm} 
\usepackage{mathtools}
\usepackage{mathrsfs}
\usepackage{siunitx}
\DeclareSIUnit{\molar}{M}

\usepackage{algorithm}
\usepackage{algpseudocode}

\usepackage{subfig}

\usepackage{tikz}

\usepackage{tabularx}
\usepackage{booktabs}

\usepackage{empheq}

\usepackage[toc,page]{appendix}

\usepackage{rotating}

\usepackage{ifthen}

\theoremstyle{plain}
\theoremstyle{plain}
\theoremstyle{plain}
\theoremstyle{plain}
\theoremstyle{plain}
\theoremstyle{definition}
\theoremstyle{remark}

\everymath{\displaystyle}
\usepackage[a4paper, total={7in, 8in}]{geometry}
\usepackage[backend=biber,style=numeric-comp,giveninits=true,doi=false,isbn=false,url=false,eprint=false,maxbibnames=99]{biblatex}
\usepackage[hidelinks]{hyperref}

\newif\ifdouble

\newcommand{\papertitle}{The role of mechano-electric feedbacks and hemodynamic coupling in scar-related ventricular tachycardia}

\newcommand{\keywordOne}{Cardiac electromechanics}
\newcommand{\keywordTwo}{Numerical simulations}
\newcommand{\keywordThree}{Mechano-electric feedback}
\newcommand{\keywordFour}{Stretch-activated channels}
\newcommand{\keywordFive}{Ventricular tachycardia}

\newcommand{\EPnumGating}{n_{\boldsymbol{w}}}

\newcommand{\ActNumVariables}{n_{\mathbf{s}}}
\newcommand{\CircNumVariables}{n_{\mathbf{c}}}

\newcommand{\PhyEP}{\mathscr{E}}

\newcommand{\PhyMec}{\mathscr{M}}

\newcommand{\PhyAct}{\mathscr{A}}

\newcommand{\PhyCirc}{\mathscr{C}}

\newcommand{\PhyCoupl}{\mathscr{V}} 

\newcommand{\Monodomain}{\mathscr{E}}

\newcommand{\VLVthreedim}{V_{\mathrm{LV}}^{\mathrm{3D}}}

\newcommand{\VLA}{V_{\mathrm{LA}}}
\newcommand{\VLV}{V_{\mathrm{LV}}}
\newcommand{\VRA}{V_{\mathrm{RA}}}
\newcommand{\VRV}{V_{\mathrm{RV}}}

\newcommand{\VnLA}{V_{\mathrm{0,LA}}}

\newcommand{\VnRA}{V_{\mathrm{0,RA}}}
\newcommand{\VnRV}{V_{\mathrm{0,RV}}}

\newcommand{\PLV}{p_{\mathrm{LV}}}

\newcommand{\PRV}{p_{\mathrm{RV}}}

\newcommand{\EpLA}{E_{\mathrm{LA}}^{\mathrm{pass}}}

\newcommand{\EpRA}{E_{\mathrm{RA}}^{\mathrm{pass}}}
\newcommand{\EpRV}{E_{\mathrm{RV}}^{\mathrm{pass}}}

\newcommand{\EaMaxLA}{E_{\mathrm{LA}}^{\mathrm{act,max}}}

\newcommand{\EaMaxRA}{E_{\mathrm{RA}}^{\mathrm{act,max}}}
\newcommand{\EaMaxRV}{E_{\mathrm{RV}}^{\mathrm{act,max}}}
\newcommand{\QarSYS}{Q_{\mathrm{AR}}^{\mathrm{SYS}}}
\newcommand{\QarPUL}{Q_{\mathrm{AR}}^{\mathrm{PUL}}}
\newcommand{\QvnSYS}{Q_{\mathrm{VEN}}^{\mathrm{SYS}}}
\newcommand{\QvnPUL}{Q_{\mathrm{VEN}}^{\mathrm{PUL}}}

\newcommand{\CarSYS}{C_{\mathrm{AR}}^{\mathrm{SYS}}}
\newcommand{\CarPUL}{C_{\mathrm{AR}}^{\mathrm{PUL}}}
\newcommand{\CvnSYS}{C_{\mathrm{VEN}}^{\mathrm{SYS}}}
\newcommand{\CvnPUL}{C_{\mathrm{VEN}}^{\mathrm{PUL}}}

\newcommand{\ParSYS}{p_{\mathrm{AR}}^{\mathrm{SYS}}}
\newcommand{\ParPUL}{p_{\mathrm{AR}}^{\mathrm{PUL}}}
\newcommand{\PvnSYS}{p_{\mathrm{VEN}}^{\mathrm{SYS}}}
\newcommand{\PvnPUL}{p_{\mathrm{VEN}}^{\mathrm{PUL}}}

\newcommand{\RarSYS}{R_{\mathrm{AR}}^{\mathrm{SYS}}}
\newcommand{\RarPUL}{R_{\mathrm{AR}}^{\mathrm{PUL}}}
\newcommand{\RvnSYS}{R_{\mathrm{VEN}}^{\mathrm{SYS}}}
\newcommand{\RvnPUL}{R_{\mathrm{VEN}}^{\mathrm{PUL}}}

\newcommand{\LarSYS}{L_{\mathrm{AR}}^{\mathrm{SYS}}}
\newcommand{\LarPUL}{L_{\mathrm{AR}}^{\mathrm{PUL}}}
\newcommand{\LvnSYS}{L_{\mathrm{VEN}}^{\mathrm{SYS}}}
\newcommand{\LvnPUL}{L_{\mathrm{VEN}}^{\mathrm{PUL}}}

\newcommand{\Rmin}{R_{\mathrm{min}}}
\newcommand{\Rmax}{R_{\mathrm{max}}}

\newcommand{\Circ}{\boldsymbol{c}}

\newcommand{\CircInit}{\boldsymbol{c}_{0}}

\newcommand{\CircRhs}{\boldsymbol{Z}}

\newcommand{\fZero}{{\mathbf{f}_0}}
\newcommand{\sZero}{{\mathbf{s}_0}}
\newcommand{\nZero}{{\mathbf{n}_0}}

\newcommand{\EPchim}{\chi_\mathrm{m}}
\newcommand{\EPCm}{C_\mathrm{m}}
\newcommand{\EPIion}{{\mathcal{I}_{\mathrm{ion}}}}
\newcommand{\EPISAC}{{\mathcal{I}_{\mathrm{SAC}}}}
\newcommand{\EPIapp}{{\mathcal{I}_{\mathrm{app}}}}

\newcommand{\EPdiffTens}{\boldsymbol{D}}

\newcommand{\EPrhsGating}{\boldsymbol{H}}

\newcommand{\EPIappReducedMax}{{\widetilde{\mathcal{I}}_{\mathrm{app}}^{\mathrm{max}}}}
\newcommand{\EPIappDuration}{t_{\mathrm{app}}}

\newcommand{\LinPlus}{{\mathrm{Lin}^+}}
\newcommand{\Inv}[1]{{\mathcal{I}_{#1}}}
\newcommand{\IIVf}{\Inv{4f}}

\newcommand{\displ}{\mathbf{d}} 

\newcommand{\BCmecCepiN}{{C_\bot^{\mathrm{epi}}}}
\newcommand{\BCmecKepiN}{{K_\bot^{\mathrm{epi}}}}
\newcommand{\BCmecCepiT}{{C_\parallel^{\mathrm{epi}}}}
\newcommand{\BCmecKepiT}{{K_\parallel^{\mathrm{epi}}}}

\newcommand{\BCmecKepiTens}{\mathbf{K}^{\mathrm{epi}}}
\newcommand{\BCmecCepiTens}{\mathbf{C}^{\mathrm{epi}}}

\newcommand{\GammaBase}{\Gamma_0^{\mathrm{base}}}
\newcommand{\GammaEpi}{\Gamma_0^{\mathrm{epi}}}
\newcommand{\GammaEndo}{\Gamma_0^{\mathrm{endo}}}

\newcommand{\GammaClosedBase}{\overline{\Gamma}_0^{\mathrm{base}}}
\newcommand{\GammaClosedEpi}{\overline{\Gamma}_0^{\mathrm{epi}}}
\newcommand{\GammaClosedEndo}{\overline{\Gamma}_0^{\mathrm{endo}}}

\newcommand{\mecF}{\mathbf{F}}
\newcommand{\mecC}{\mathbf{C}}
\newcommand{\tenspiola}{\mathbf{P}}
\newcommand{\identity}{\mathbf{I}}

\newcommand{\mecNref}{{\mathbf{N}}}

\newcommand{\BCmecVbaseLV}{\mathbf{v}_{\mathrm{LV}}^{\mathrm{base}}}

\newcommand{\CaiFactor}{\omega_{\mathrm{Ca}}}

\newcommand{\Cai}{{[\mathrm{Ca}^{2+}]_{\mathrm{i}}}}

\newcommand{\SL}{{SL}} 

\newcommand{\Tens}{T_\mathrm{a}} 
\newcommand{\ActStateHF}{\mathbf{s}}

\newcommand{\ActRhs}{\boldsymbol{K}}

\newcommand{\Koff}{{k_\text{off}}}
\newcommand{\Kbasic}{{k_\text{basic}}}
\newcommand{\Kd}{{k_\text{d}}}
\newcommand{\KdZero}{{\overline{k}_\text{d}}}
\newcommand{\KdAlpha}{{\alpha_\Kd}}

\newcommand{\aXB}{{a_{\text{XB}}}}
\newcommand{\tmP}{\mathcal{P}}
\newcommand{\fP}{{f_\tmP}}

\newcommand{\OmegaRef}{\widetilde{\Omega}}

\newcommand{\pressLVRef}{\widetilde{p}_{\mathrm{LV}}}

\newcommand{\tensRef}{\widetilde{T}_{\mathrm{a}}}

\newcommand{\Gs}{G_{\text{s}}}
\newcommand{\urev}{u_{\text{rev}}}
\newcommand{\Vrev}{V_{\text{rev}}}

\graphicspath{{pictures/}}
\doublefalse

\title{{\papertitle}}
\author{Matteo Salvador$^1$,
	      Francesco Regazzoni$^1$,
		    Stefano Pagani$^1$,\\
		    Luca Dede'$^1$,
				Natalia Trayanova$^2$,
		    Alfio Quarteroni$^{1, 3}$}

\date{\footnotesize
      $^1$ MOX-Dipartimento di Matematica, Politecnico di Milano, Milan, Italy \\
	    $^2$ Department of Biomedical Engineering, Johns Hopkins University, Baltimore, MD, USA\\
	    $^3$ Professor Emeritus, \'Ecole Polytechnique F\'ed\'erale de Lausanne, Lausanne, Switzerland\\[2ex]}

\begin{document}
    \newif\ifCIBM
    \CIBMfalse

	\maketitle

	\begin{abstract}
		Mechano-electric feedbacks (MEFs), which model how mechanical stimuli are transduced into electrical signals, have received sparse investigation by considering electromechanical simulations in simplified scenarios.
In this paper, we study the effects of different MEFs modeling choices for myocardial deformation and nonselective stretch-activated channels (SACs) in the monodomain equation.
We perform numerical simulations during ventricular tachycardia (VT) by employing a biophysically detailed and anatomically accurate 3D electromechanical model for the left ventricle (LV) coupled with a 0D closed-loop model of the cardiocirculatory system.
We model the electromechanical substrate responsible for scar-related VT with a distribution of infarct and peri-infarct zones.
Our mathematical framework takes into account the hemodynamic effects of VT due to myocardial impairment and allows for the classification of their hemodynamic nature, which can be either stable or unstable.
By combining electrophysiological, mechanical and hemodynamic models, we observe that all MEFs may alter the propagation of the action potential and the morphology of the VT.
In particular, we notice that the presence of myocardial deformation in the monodomain equation may change the VT basis cycle length and the conduction velocity but do not affect the hemodynamic nature of the VT.
Finally, nonselective SACs may affect wavefront stability, by possibly turning a hemodynamically stable VT into a hemodynamically unstable one and vice versa.

	\end{abstract}


	\section{Introduction}
\label{sec: introduction}

Cardiac arrhythmias result from an irregular electrical activity of the human heart.
Among them, ventricular tachycardia (VT), which manifests with an accelerated heart rate, is one of the most life-threatening rhythm disorders. 
VT may be classified as either hemodynamically stable or unstable, depending on the capability of the heart to effectively pump blood in the circulatory system.
In the former case antiarrhythmic drugs are generally employed, while in the latter case cardioversion is 
needed \cite{Epstein2008}.
According to the specific pathogenesis, the stability of the VT remains the same or changes over time.
Moreover, it may also degenerate towards ventricular fibrillation (VF), a life-threatening condition in which the ventricular activity is fully disorganized and chaotic, leading to heart failure \cite{Samie2001}.

In the clinical framework, these pathological scenarios can be hardly ever fully investigated and predicted for all patients.
For this reason, biophysically detailed computational heart models could be used to provide a deeper understanding of the hemodynamic response to VT and 
to characterize the electromechanical substrate leading to dangerous arrhythmias.

While electrophysiological simulations are well-established for scar-related VT identification and treatment on human ventricles \cite{Arevalo2016, Cedilnik2018, Deng2019, Prakosa2018}, patient-specific electromechanical models coupled with closed-loop cardiovascular circulation have been just recently used for the first time to enhance our knowledge on VT \cite{Salvador2021}.
Indeed, on one side, the physiological processes by which the mechanical behavior alters the electrical activity of the human heart, known as mechano-electric feedbacks (MEFs), are relevant and not fully elucidated \cite{Kohl2004, Kohl2013, ColliFranzone2017, Timmermann2017AnIA, Taggart1999, Panfilov2010, Varela2021}.
On the other hand, the identification of the hemodynamic nature of the VT has significant clinical implications \cite{Salvador2021}.


The goal of this study is threefold:
(1) to combine electrophysiology, activation, mechanics and hemodynamics in several numerical simulations of scar-related VT, with the aim of uncovering the roles of myocardial deformation and the recruitment of SACs on VT stability;
(2) to show that our computational model effectively reproduces both hemodynamically stable and hemodynamically unstable VT;
(3) to capture relevant microscopic mechanisms during VT, such as the incomplete relaxation of sarcomeres, given the biophysical detail of our electromechanical model.

	\section{Mathematical models}
\label{sec: mathematicalmodeling}
We provide an overview of the 3D cardiac electromechanical model coupled with a 0D closed-loop model of the cardiovascular system (see \cite{regazzoni2020model, Salvador2021}).

\begin{figure}[t!]
	\centering
	\includegraphics[keepaspectratio, width=0.45\textwidth]{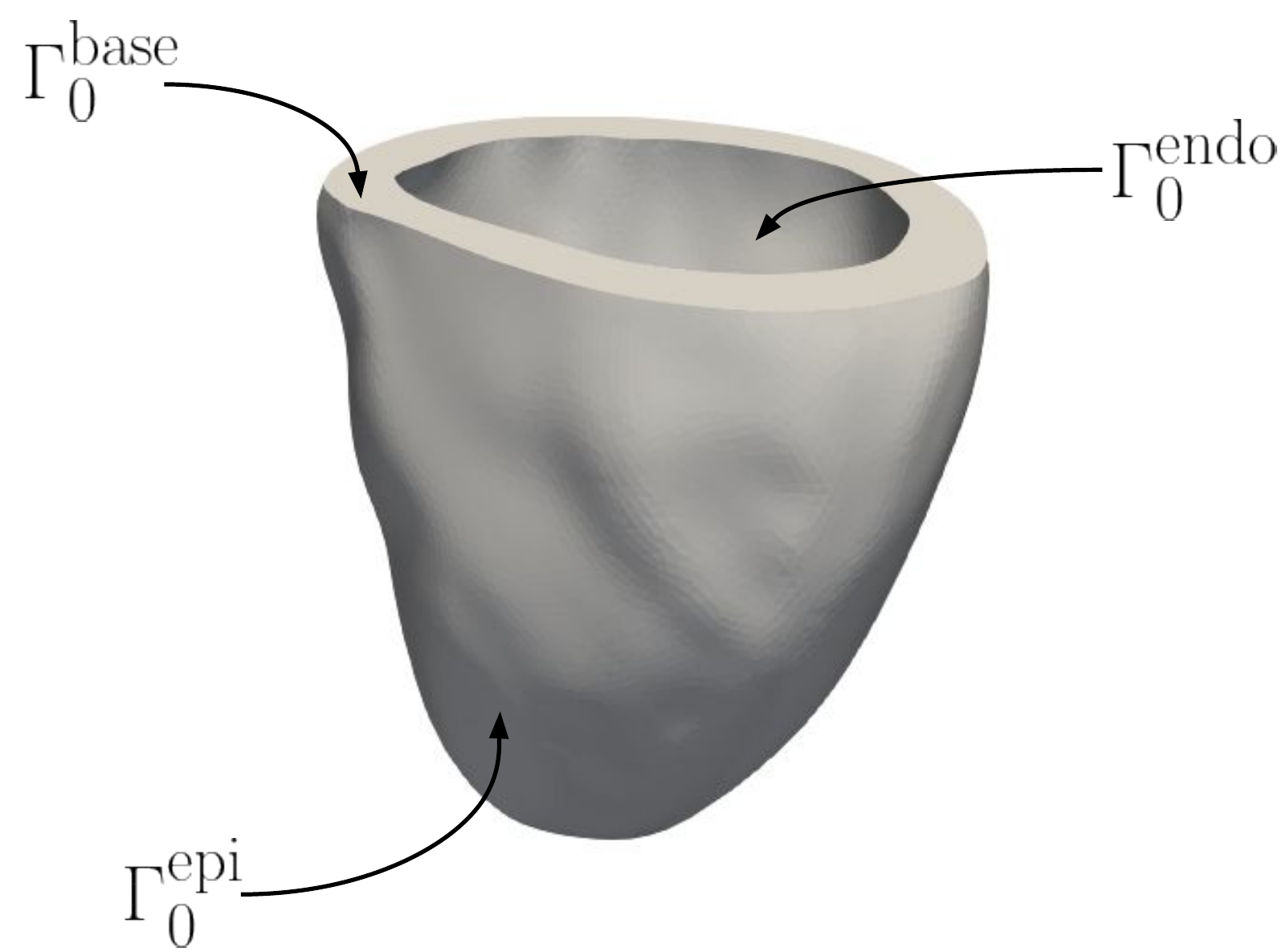}
	\caption{Computational domain $\Omega_0$, given by a human LV.}
	\label{fig: tags_LV}
\end{figure}

We depict in Figure~\ref{fig: tags_LV} the computational domain $\Omega_0$, which represents a human left ventricle (LV) taken from the Zygote Solid 3D heart model \cite{Zygote}.
Its boundary $\partial \Omega_0$ is partitioned into epicardium $\GammaEpi$, endocardium $\GammaEndo$ and base $\GammaBase$. We recall that $\partial \Omega_0=\GammaClosedEpi \cup \GammaClosedEndo \cup \GammaClosedBase$.

\subsection{3D-0D closed-loop electromechanical model}
\label{subsec: elettromech_model}

We consider a multiphysics and multiscale 3D-0D mathematical framework comprised of four different core models, namely cardiac electrophysiology $(\PhyEP)$ \cite{franzone2014mathematical,colli2018numerical,nobile2012active,levrero2020sensitivity}, mechanical activation $(\PhyAct)$ \cite{niederer2011length,rossi2014thermodynamically,ruiz2014mathematical,land2017model,regazzoni2020biophysically}, passive mechanics $(\PhyMec)$ \cite{guccione1991passive,guccione1991finite,ogden1997non,holzapfel2009constitutive} and blood circulation $(\PhyCirc)$ \cite{quarteroni2016geometric,blanco20103d,kerckhoffs2007coupling,hirschvogel2017monolithic,augustin2021,arts2005adaptation,regazzoni2020model}.
The volume conservation constraint~$(\PhyCoupl)$ defines the coupling condition between the LV cardiac electromechanics and the remaining part of cardiovascular system.~\cite{regazzoni2020model}.
All these models represent several physiological processes, ranging from the cellular level to the organ scale.

Given the computational domain $\Omega_0$ and the time interval $t \in (0,T]$, the mathematical model reads:
\input{EMmodel_lv.tex} in which the model unknowns are:
\begin{equation} \label{eqn: unknonws}
\begin{aligned}
	u \colon &\Omega_0 \times (0,T] \to \mathbb{R} , &
	\boldsymbol{w} \colon &\Omega_0 \times (0,T] \to \mathbb{R}^{\EPnumGating} , \\
	\ActStateHF \colon &\Omega_0 \times (0,T] \to \mathbb{R}^{\ActNumVariables} , &
	\displ \colon &\Omega_0 \times (0,T] \to \mathbb{R}^3, &
	\Circ \colon &(0,T] \to \mathbb{R}^{\CircNumVariables} , \\
	\PLV \colon &(0,T] \to \mathbb{R},
\end{aligned}
\end{equation}
\noindent
where $u$ is the transmembrane potential, $\boldsymbol{w}$ is the vector containing both gating and ionic variables, $\ActStateHF$ are the states variables of the active force generation model, $\displ$ represents the tissue mechanical displacement, $\Circ$ defines the state vector of the circulation model (including pressures, volumes and fluxes of the different compartments of the cardiocirculatory system) and $\PLV$ and $\PRV$ are the left and right ventricular pressure, respectively.
Finally, all variables are endowed with suitable initial conditions in $\Omega_0 \times \{0\}$:

\begin{equation} \label{eqn: initial_conditions}
u = u_0, \quad \boldsymbol{w} = \boldsymbol{w}_0, \quad \ActStateHF = \ActStateHF_0, \quad \displ = \displ_0, \quad \dfrac{\partial \displ}{\partial t} = \Dot{\displ}_0, \quad \Circ=\CircInit.
\end{equation}


\subsubsection{Electrophysiology $(\PhyEP)$}
\label{subsec: electrophysiology}
We model the electrical activity of the myocardium by means of Eq.~\eqref{eqn: EP}, that is the monodomain equation coupled with a suitable ionic model for the human ventricular action potential ~\cite{BuenoOrovio2008, franzone2014mathematical, LuoRudy1994, ten2006alternans}.

In the electrophysiological model $(\PhyEP)$, $\EPchim$ is the surface area-to-volume ratio of cardiomyocytes, $\EPCm$ represents the transmembrane capacitance per unit area.
The applied current $\EPIapp$ mimics the effect of the Purkinje network \cite{vergara2014patient,costabal2016generating,landajuela2018numerical} by triggering the action potential at specific locations of the myocardium.
The reaction terms $\EPIion$ and $\EPrhsGating$ (specified by the ionic model at hand) couple the action potential propagation to the cellular dynamics. Specifically, we use the ten Tusscher-Panfilov ionic model (TTP06), which is able to accurately describe ions dynamics across the cell membrane \cite{ten2006alternans}.
Furthermore, Eq.~\eqref{eqn: EP} is equipped with homogeneous Neumann boundary conditions for $u$ at the boundary $\partial \Omega_0$, which defines the condition of electrically isolated domain for $\Omega_0$.
$\EPdiffTens$ represents the conductivity tensor, that reads:
\begin{equation} \label{eqn: diff_tensor}
\EPdiffTens = \eta \sigma_{\text{l}} \frac{\mecF\fZero \otimes \mecF\fZero}{\|\mecF \fZero \|^2} + \eta \sigma_{\text{t}} \frac{\mecF\sZero \otimes \mecF\sZero}{\|\mecF \sZero \|^2} + \eta \sigma_{\text{n}} \frac{\mecF\nZero \otimes \mecF\nZero}{\|\mecF \nZero \|^2},
\end{equation}
where $\sigma_{\text{l}}, \sigma_{\text{t}}, \sigma_{\text{n}}$ are the longitudinal, transversal and normal conductivities, respectively.
The parameter $\eta = \eta(\boldsymbol{x}) \in [0, 1]$ takes into account the effect of scars, grey zones and non-remodeled regions. This parameter is also incorporated inside the formulation of the TTP06 model, as described in \cite{Salvador2021}.

The deformation tensor $\mecF = \identity + \nabla \displ$ and its determinant $J = \det(\mecF)$ are needed to perform the mechano-electric coupling \cite{Trayanova2011}.
Indeed, we model the so called mechano-electric feedbacks (MEFs) \cite{Costabal2017, Hazim2021}.
The geometry-mediated MEFs incorporate the effects of displacement $\displ$ on the cardiac tissue, while other physiological processes act at the level of single cardiomyocytes \cite{Kohl2004, levrero2020sensitivity}.
Among them, some examples are selective (e.g. $K^\text{+}$-permeable) or nonselective SACs and intracellular calcium $\Cai$ binding to sarcolemmal buffers, the latter requiring more sophisticated ventricular ionic models \cite{Bartolucci2020, Tomek2019}.
In this paper we model nonselective SACs by means of the following formulation \cite{levrero2020sensitivity}:
\begin{equation} \label{eqn: SAC}
\EPISAC(u, \mecF) = \Gs (\|\mecF \fZero \| - 1)_{+} (u - \urev),
\end{equation}
where $\Gs$ and $\urev$ represent the conductance of the channels and the reversal potential, respectively.
SACs may alter the shape of the action potential (AP) by lengthening or shortening its duration (APD) and by generating higher resting potentials.
This may induce early or delayed afterdepolarizations (EADs or DADs) and premature excitation.

Both geometry-mediated MEFs and nonselective SACs are generally known to be pro-arrhythmic in pathological scenarios, because they increase the likelihood of having extra stimuli during cardiac contraction and relaxation \cite{Kohl2004}.

In this work we consider Eq.~\eqref{eqn: EP} with several degrees of complexity to assess similarities and differences in the outcomes of the electromechanical simulations during VT.
We report in Tab.~\ref{tab: EP_models} the mathematical models and the parametrizations that we consider.
In particular, the choice of three different mathematical models for the geometry-mediated MEFs, namely $(\Monodomain_{\text{gMEF-minimal}})$, $(\Monodomain_{\text{gMEF-enhanced}})$ and $(\Monodomain_{\text{gMEF-full}})$, is motivated by the numerous formulations of this type of feedback that can be found in the literature \cite{Chapelle2009, Collin2019, levrero2020sensitivity, quarteroni2017integrated}.
Indeed, we range from minimal to complete inclusion of geometry-mediated MEFs in the monodomain equation.

\begin{table}[t!]
 \centering
 \caption{Modeling choices for the monodomain equation that have been used in this paper. We consider different parametrizations for $\EPISAC(u, \mecF)$ in terms of $\Gs$ and $\urev$. $\EPdiffTens_{\mathbf{I}}$ indicates the conductivity tensor in Eq.~\eqref{eqn: diff_tensor} with $\mecF = \identity$.}
 \label{tab: EP_models}
 \renewcommand{\arraystretch}{2.5}
 \begin{tabular}{ll}
 \hline\noalign{\smallskip}
 Model name & Equation \\
 \noalign{\smallskip}\hline\noalign{\smallskip}
 $(\Monodomain)$ & \( \EPchim \left[ \EPCm \dfrac{\partial u}{\partial t} + \EPIion(u, \boldsymbol{w}) \right] - \nabla \cdot ( \EPdiffTens_{\mathbf{I}} \nabla u) = \EPchim \EPIapp(t) \) \\
 $(\Monodomain_{\text{gMEF-minimal}})$ & \( \EPchim \left[ \EPCm \dfrac{\partial u}{\partial t} + \EPIion(u, \boldsymbol{w}) \right] - \nabla \cdot ( J \mecF^{-1} \EPdiffTens_{\mathbf{I}} \, \mecF^{-T} \nabla u) = \EPchim \EPIapp(t) \) \\
 $(\Monodomain_{\text{gMEF-enhanced}})$ & \( \EPchim \left[ \EPCm \dfrac{\partial u}{\partial t} + \EPIion(u, \boldsymbol{w}) \right] - \nabla \cdot ( J \mecF^{-1} \EPdiffTens \, \mecF^{-T} \nabla u) = \EPchim \EPIapp(t) \) \\
 $(\Monodomain_{\text{gMEF-full}})$ & \( J \EPchim \left[ \EPCm \dfrac{\partial u}{\partial t} + \EPIion(u, \boldsymbol{w}) \right] - \nabla \cdot ( J \mecF^{-1} \EPdiffTens \, \mecF^{-T} \nabla u) = J \EPchim \EPIapp(t) \) \\
 $(\Monodomain_{\text{SAC}})$ & \( \EPchim \left[ \EPCm \dfrac{\partial u}{\partial t} + \EPIion(u, \boldsymbol{w}) + \EPISAC(u, \mecF) \right] - \nabla \cdot ( \EPdiffTens_{\mathbf{I}} \nabla u) = \EPchim \EPIapp(t) \) \\
 $(\Monodomain_{\text{gMEF-full, SAC}})$ & \( J \EPchim \left[ \EPCm \dfrac{\partial u}{\partial t} + \EPIion(u, \boldsymbol{w}) + \EPISAC(u, \mecF) \right] - \nabla \cdot ( J \mecF^{-1} \EPdiffTens \, \mecF^{-T} \nabla u) = J \EPchim \EPIapp(t) \) \\
 \noalign{\smallskip}\hline
 \end{tabular}
\end{table}

\subsubsection{Activation $(\PhyAct)$}
\label{subsec: mechanicalactivation}

To model how the calcium wave following the AP triggers a series of chemo-mechanical reactions within sarcomeres, resulting in the generation of an active force in the muscle, we use the mean-field version of the model proposed in \cite{regazzoni2020biophysically}, henceforth denoted by RDQ20-MF.
This mathematical model is based on a biophysically accurate description of regulatory and contractile proteins and their dynamics.
Thanks to suitable dimensionality reduction techniques, the dynamics of stochastic processes underlying both chemical and mechanical microscale transitions are described in only 20 ODEs (that is $\ActStateHF(t) \in \mathbb{R}^{20}$).
The computational cost of this model is thus comparable to that of phenomenological models \cite{niederer_hunter_smith_2006,rice2008approximate,land2017model}, while providing a biophysically detailed description that is consistent with the level of detail and mechanistic understanding of the ionic model needed for the present study.

The inputs of the RDQ20-MF model are the intracellular calcium concentration $\Cai$ coming from $(\PhyEP)$, the sarcomere length (denoted by $\SL$) and its time derivative (the latter allows to account for the so-called force-velocity relationship \cite{keener2009mathematical}).
The variable $\SL$ is obtained as $\SL=\SL_0 \sqrt{\IIVf}$, where $\SL_0$ is the sarcomere slack length and the fourth invariant $\IIVf = \mecF \fZero \cdot \mecF \fZero$ measures the tissue stretch in the fiber direction.
Finally, the output of the RDQ20-MF model is the active tension generated at the microscale, that can be obtained as $\Tens(\ActStateHF,\SL)$.

Differently from \cite{Salvador2021}, here we consider a biophysically detailed active stress model.
For this reason, we do not need to put the parameter $\eta$ inside its parametrization.
Indeed, the active force generation mechanisms are properly handled by the active stress model, which receives different intracellular calcium waves from $(\PhyEP)$ according to the specific area of the myocardium (scar, grey zone or healthy) and provides physiological values of active tension in all cases.

\subsubsection{Mechanics $(\PhyMec)$}
\label{subsec: actpassmechanics}
We employ the momentum conservation equation reported in Eq.~\eqref{eqn: Mec} to model the dynamics of the displacement $\displ$ of the myocardium. We also consider the hyperelasticity assumption, so that the strain energy function can be differentiated with respect to the deformation tensor $\mecF$ to obtain $\tenspiola$ \cite{ogden1997non,guccione1991finite}.

In the mechanical model $(\PhyMec)$, $\rho_{\text{s}}$ represents the density of the myocardium.
The Piola-Kirchhoff stress tensor $\tenspiola = \tenspiola(\displ, \Tens)$ is additively decomposed according to:
\begin{equation} \label{eqn: piola}
\displaystyle \tenspiola(\displ, \Tens)
= \dfrac{\partial \mathcal{W}(\mecF)}{\partial \mecF}
+ \Tens(\ActStateHF, \SL) \frac{\mecF \fZero \otimes \fZero}{\sqrt{\IIVf}},
\end{equation}
where the first term stands as the passive part of the tensor $\tenspiola$, while the latter as the active one; $\mathcal{W} : \LinPlus \to \mathbb{R}$ is the strain energy density function, $\Tens(\ActStateHF, \SL)$ is the active tension, provided by the activation model $(\PhyAct)$.

To model the passive behaviour of cardiac tissue, we employ the orthotropic Guccione constitutive law \cite{guccione1991passive}, whose strain energy function is defined as:
\begin{equation}
\displaystyle \mathcal{W} = \tfrac{\kappa}{2} \left( J - 1 \right) \text{log}(J) + \tfrac{a}{2} \left( e^Q  - 1 \right),
\end{equation}
where the first term is the volumetric energy with the bulk modulus $\kappa$, which penalizes large variation of volume to enforce a weakly incompressible behavior \cite{peng1997compressible,doll2000development}, and the latter is the deviatoric energy where $a$ is the stiffness scaling parameter and the exponent $Q$ reads:
\begin{equation}
\begin{split}
Q &= b_{\text{ff}} E_{\text{ff}}^2  + b_{\text{ss}} E_{\text{ss}}^2 + b_{\text{nn}} E_{\text{nn}}^2+ b_{\text{fs}} \left( E_{\text{fs}}^2 + E_{\text{sf}}^2 \right) + b_{\text{fn}} \left( E_{\text{fn}}^2 + E_{\text{nf}}^2 \right) + b_{\text{sn}} \left( E_{\text{sn}}^2 + E_{\text{ns}}^2 \right),
\end{split}
\end{equation}
in which $E_\text{ij} = \boldsymbol{E} \boldsymbol{\mathrm{i}}_\text{0} \cdot \boldsymbol{\mathrm{j}}_\text{0} \;\;\; \text{for} \;\;\; \mathrm{i, j} \in \{ \mathrm{f, s, n} \}$
are the entries of $\textbf{E} = \tfrac{1}{2} \left( \mecC - \identity \right)$, i.e the Green-Lagrange strain energy tensor, being $\mecC = \mecF^{T} \mecF$ the right Cauchy-Green deformation tensor.

To model the mechanical effects of the pericardial sac \cite{Gerbi2018monolithic,pfaller2019importance,strocchi2020simulating}, we impose at the epicardial boundary $\GammaEpi$ a generalized Robin boundary condition $\tenspiola(\displ, \Tens) \mecNref = \BCmecKepiTens \displ + \BCmecCepiTens \tfrac{\partial \displ}{\partial t}$ by defining the tensors $		\BCmecKepiTens =
\BCmecKepiT (\mecNref \otimes \mecNref - \identity) -
\BCmecKepiN (\mecNref \otimes \mecNref)$ and $		\BCmecCepiTens =
\BCmecCepiT (\mecNref \otimes \mecNref - \identity) -
\BCmecCepiN (\mecNref \otimes \mecNref),$
where $\BCmecKepiN$, $\BCmecKepiT$, $\BCmecCepiN$, $\BCmecCepiT \in \mathbb{R}^+$ are the stiffness and viscosity parameters of the epicardial tissue in the normal and tangential directions, respectively. Normal stress boundary conditions were imposed at the endocardium of the LV $\GammaEndo$, in which $\PLV(t)$ is the pressure exerted by the blood in the LV, modeled by means of the 0D closed-loop circulation model \cite{regazzoni2020model}.
To take into account the effect of the neglected part, over the basal plane, on the ventricular domain, we set on $\GammaBase$ the energy consistent boundary condition $\tenspiola(\displ, \Tens) \mecNref = \displaystyle | J \mecF^{-T} \mecNref | \BCmecVbaseLV(t)$ \cite{regazzoni2020machine}, where:
\begin{equation}
\BCmecVbaseLV(t) = \displaystyle \frac{\int_{\GammaEndo} \PLV(t) J \mecF^{-T} \mecNref d \Gamma_0}{\int_{\GammaBase} | J \mecF^{-T} \mecNref | d \Gamma_0}.
\end{equation}

\subsubsection{Blood circulation $(\PhyCirc)$ and coupling conditions $(\PhyCoupl)$ }
\label{sec: mathematicalmodeling: circulation}
We model the blood circulation through the entire cardiovascular system by means of a closed-loop model recently proposed in \cite{regazzoni2020model}.
In the 0D closed-loop model, systemic and pulmonary circulations are modeled with RLC circuits, heart chambers are described by time-varying elastance elements and non-ideal diodes stand for the heart valves \cite{regazzoni2020model}.

We represent the circulation core model $(\PhyCirc)$ with a system of ODEs, which is expressed by Eq.~\eqref{eqn: Circ}, where $\CircRhs$ is a proper function (defined in \cite{regazzoni2020model}) and $\Circ(t)$ includes pressures, volumes and fluxes of the different compartments composing the vascular network:
\begin{equation*}
	\begin{split}
		\Circ(t) =
		(&\VLA(t), \VLV(t), \VRA(t), \VRV(t), \ParSYS(t), \PvnSYS(t), \ParPUL(t), \PvnPUL(t), \\&\QarSYS(t), \QvnSYS(t), \QarPUL(t), \QvnPUL(t))^T.
	\end{split}
\end{equation*}
Here $\VLA$, $\VRA$, $\VLV$ and $\VRV$ refer to the volumes of left atrium, right atrium, left ventricle and right ventricle, respectively;
$\ParSYS$, $\QarSYS$, $\PvnSYS$, $\QvnSYS$, $\ParPUL$, $\QarPUL$, $\PvnPUL$ and $\QvnPUL$ express pressures and flow rates of the systemic and pulmonary circulation (arterial and venous).
For the complete mathematical description of the 0D circulation lumped model we refer to \cite{regazzoni2020model}.
To couple the 0D circulation model $(\PhyCirc)$ with the 3D LV model, given by $(\PhyEP)$--$(\PhyAct)$--$(\PhyMec)$, we follow the same strategy proposed in \cite{regazzoni2020model}.
In particular, we replace the time-varying elastance elements representing the LV in $(\PhyCirc)$ with its corresponding 3D electromechanical description and we introduce the coupling condition $(\PhyCoupl)$ where:
\begin{equation*}
	\displaystyle
	\VLVthreedim(\displ(t)) =
	\int_{\GammaEndo} J(t)
	\left(\left( \mathbf{h} \otimes \mathbf{h}  \right) \left(\mathbf{x} + \displ(t) - \mathbf{b} \right) \right)
	\cdot \mecF^{-T}(t) \mecNref \, d \Gamma_0
\end{equation*}
in which $\mathbf{h}$ is a vector orthogonal to the LV centerline (i.e. lying on the ventricular base) and $\mathbf{b}$ lays inside the LV \cite{regazzoni2020model}.

Due to $(\PhyCoupl)$, in the 3D-0D coupled model~\eqref{eqn: EM}, $\PLV(t)$ is not determined by the 0D circulation model~\eqref{eqn: Circ}, but rather acts as Lagrange multipliers and enforces the constraint $(\PhyCoupl)$.

\subsection{Reference configuration}
\label{sec: refconf}
Cardiac geometries are acquired from in vivo medical images through imaging techniques.
These geometries are in principle not stress free, mainly because there is always a pressure acting on the endocardium.
Therefore, one needs to estimate the unloaded (i.e. stress-free) configuration (also named reference configuration) to which the 3D-0D cardiac electromechanical model~\eqref{eqn: EM} refers. To recover the reference configuration $\Omega_0$, starting from a geometry acquired from medical images $\OmegaRef$, we apply the same procedure proposed for the left ventricle in \cite{regazzoni2020model}.

Specifically, we solve an inverse problem: 
find the solution $\displ = \displ_{\Omega_0}$ of the following differential problem
\begin{equation}
\begin{cases}
	\nabla \cdot \tenspiola(\displ, \Tens) = \boldsymbol{0} & $in$ \; \Omega_0 , \\
	\tenspiola(\displ, \Tens) \mecNref +
	\BCmecKepiTens \displ
	= \mathbf{0}
	& $on$ \; \GammaEpi , \\
\tenspiola(\displ, \Tens) \mecNref = -\PLV(t) \, J \mecF^{-T} \mecNref & $on$ \; \GammaEndo, \\
\tenspiola(\displ, \Tens) \mecNref = \displaystyle | J \mecF^{-T} \mecNref | \BCmecVbaseLV(t)
& $on$ \; \GammaBase,
\end{cases}
\label{eqn: mechanics_steadystate}
\end{equation}
obtained for $\PLV = \pressLVRef$ and $\Tens = \tensRef$, such that we get the domain $\OmegaRef = \{ \mathbf{x} : \mathbf{x} = \mathbf{x}_\text{0} + \displ \;\; \forall \mathbf{x}_\text{0} \in \Omega_0 \}$.

Finally, to properly initialize the numerical simulation, we inflate the ventricular reference configuration $\Omega_0$ by solving problem~\eqref{eqn: mechanics_steadystate}, where we set the pressures $\PLV=p_\text{LV}^\text{ED}$.
The value $p_\text{LV}^\text{ED}$ is chosen to bring the left ventricle to a defined volume $V_\text{LV}^\text{ED}$, as explained in \cite{Piersanti2021EM}.
Then, the solution $\displ$ of the problem~\eqref{eqn: mechanics_steadystate} is set as initial condition $\displ_0$ for $\displ$ (with $\Dot{\displ}_0 = \mathbf{0}$) in $(\PhyMec)$.

	\section{Numerical methods}
\label{sec: methods}

\begin{figure}[t!]
\centering
\includegraphics[keepaspectratio, width=0.6\textwidth]{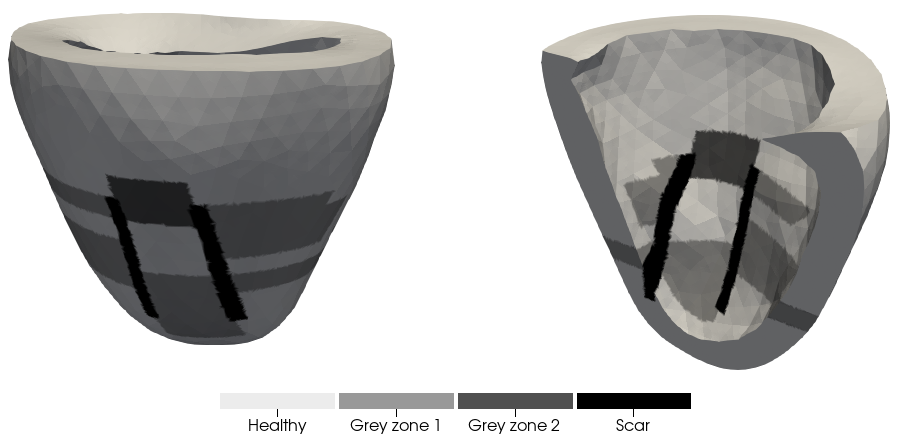}
\caption{Zygote LV with the idealized distribution of scars (black), grey zones (grey) and non-remodeled regions (white) over the myocardium. Volumetric view (left) and cut view (right). The first type of grey zone corresponds to $\eta = 0.2$, while on the second one $\eta = 0.1$ is prescribed.}
\label{fig: ischemicregion}
\end{figure}

We depict the geometric model of the Zygote LV in Fig.~\ref{fig: ischemicregion}.
We consider an idealized distribution of ischemic regions, which is made by two scars and two different types of grey zones.

To generate the fibers distribution (field $\fZero$) for our geometry, we use the Bayer-Blake-Plank-Trayanova algorithm \cite{bayer2012novel, piersanti2021modeling} with $\alpha_{\text{epi}} = -60^{\circ}$, $\alpha_{\text{endo}} = 60^{\circ}$, $\beta_{\text{epi}} = 20^{\circ}$, and $\beta_{\text{endo}} = -20^{\circ}$.

We employ the strategy proposed in \cite{Regazzoni2021} to choose proper initial conditions for the 3D electromechanical simulations, by relying on a 0D emulator.
In particular, to reach a steady state for the electrophysiological variables, we trigger periodic stimuli with a period equal to $\SI{0.45}{\second}$, for a total duration of $\SI{450}{\second}$.
On the other hand, for the 3D numerical simulations, we apply an S1-S2 stimulation protocol consisting of two gaussian stimuli, the first one applied at time $t = \SI{0}{\second}$ and the second one applied at time $t = \SI{0.45}{\second}$.
In this way the 3D electromechanical simulations are a natural continuation of the 0D initialization, and the VT can be directly induced by just using a single S1 before the S2.

We use a segregated-intergrid-staggered scheme to numerically discretize the electromechanical model \cite{Piersanti2021EM}.
Indeed, we solve $(\PhyEP)$, $(\PhyAct)$, $(\PhyMec)-(\PhyCoupl)$ and $(\PhyCirc)$ sequentially, by employing different time and space resolutions according to the specific core model.

We employ the Finite Element Method (FEM) for the space discretization of electrophysiological, activation and mechanical models \cite{quarteroni2009numerical, quarteroni2019cardiovascular}. We use hexahedral meshes and $\mathbb{Q}_1$ finite element spaces for the space discretization of all the core models.
This multiphysics problem presents different space resolutions according to the specific model at hand \cite{regazzoni2020model, regazzoni2020numerical}.
In the framework of cardiac electrophysiology, we consider a fine geometrical description (1'987'285 DOFs, 1'926'912 cells, $h_{\text{mean}} \approx \SI{0.86}{\milli\meter}$) to accurately capture the electric propagation due to fine-scale phenomena arising from the continuum modeling of the cellular level, especially with the aim of reproducing and properly address arrhythmias \cite{salvador2020intergrid}.
On the other hand, cardiac mechanics allows a lower space resolution (35'725 DOFs, 30'108 cells, $h_{\text{mean}} \approx \SI{3.3}{\milli\meter}$). This eases its numerical solution, which is computationally demanding, especially for the assembling phase \cite{salvador2020intergrid}.
Differently from \cite{Salvador2021}, in this paper we consider a significant space scales separation between electrophysiology and mechanics because we are dealing with a simple distribution of ischemic regions.

For time discretization, we use backward differentiation formula (BDF) schemes \cite{quarteroni2009numerical}. In particular, we employ a third order BDF scheme for electrophysiology and a first order BDF scheme for activation, mechanics and circulation.
These choices allow to accurately capture the fast time dynamics of the model variables without unbearable restrictions on the time step, while not introducing numerical instabilities.
We treat the nonlinear terms coming from both electrophysiological and activation models in a semi-implicit fashion.
Cardiac mechanics is numerically advanced in time with a fully implicit scheme.
The cardiocirculatory model is solved with an explicit method \cite{Piersanti2021EM}.
To stabilize the numerical oscillations arising from the feedback of the tissue mechanics on the activation model, without having to resort to a monolithic strategy, we rely on the scheme proposed in \cite{Regazzoni2021oscillation}.
We use a smaller time step for electrophysiology ($\tau = \SI{50}{\mu \second}$) than for activation, mechanics and circulation ($\Delta t = \SI{500}{\mu \second}$).
We set the final time $T = \SI{4}{\second}$ for all the numerical simulations.


The mathematical models of Sec.~\ref{sec: mathematicalmodeling} and the numerical methods presented in this section have been implemented in \texttt{life\textsuperscript{x}} (\url{https://lifex.gitlab.io/lifex}), a high-performance \texttt{C++} library developed within the iHEART project and based on the \texttt{deal.II} (\url{https://www.dealii.org}) Finite Element core \cite{dealII92}.

The numerical simulations were performed on a HPC facility available at MOX for the iHEART project.
The entire cluster is endowed with 8 Intel Xeon Platinum 8160 processors, for a total of 192 computational cores and a total amount of 1.5TB of available RAM.

	\section{Numerical results}
\label{sec: numericalresults}

We present electromechanical simulations to evaluate the effects of MEFs.
We consider different modeling choices for the monodomain equation, as reported in Tab.~\ref{tab: EP_models}.
We start from a baseline simulation with model $(\Monodomain)$, in which we obtain a stable VT.
Then, we compare the effects of different geometry-mediated MEFs, i.e. models $(\Monodomain)$, $(\Monodomain_{\text{gMEF-minimal}})$, $(\Monodomain_{\text{gMEF-enhanced}})$ and $(\Monodomain_{\text{gMEF-full}})$.
We also study the impact of different parametrizations for SACs, i.e. $(\Monodomain_{\text{SAC}})$, with respect to $(\Monodomain)$.
Finally, we evaluate the combined effects of geometry-mediated MEFs and nonselective SACs.
We report in Appendix~\ref{app:params} the values of the parameters that we use to get the numerical results that will be discussed in this paper.

\subsection{Baseline simulation}
\label{subsec: baseline}

\newcommand{\EMbaselineVT}[2]{
	\subfloat[][$t = \SI{#2}{\second}$]{\includegraphics[width=0.45\textwidth]{pictures/EM_VT_none/EM_VT_none_#1.png}}}

\ifCIBM \begin{figure}[h!] \else \begin{figure}[p] \fi
	\center
	\captionsetup[subfigure]{labelformat=empty}
	\includegraphics[width=0.45\textwidth]{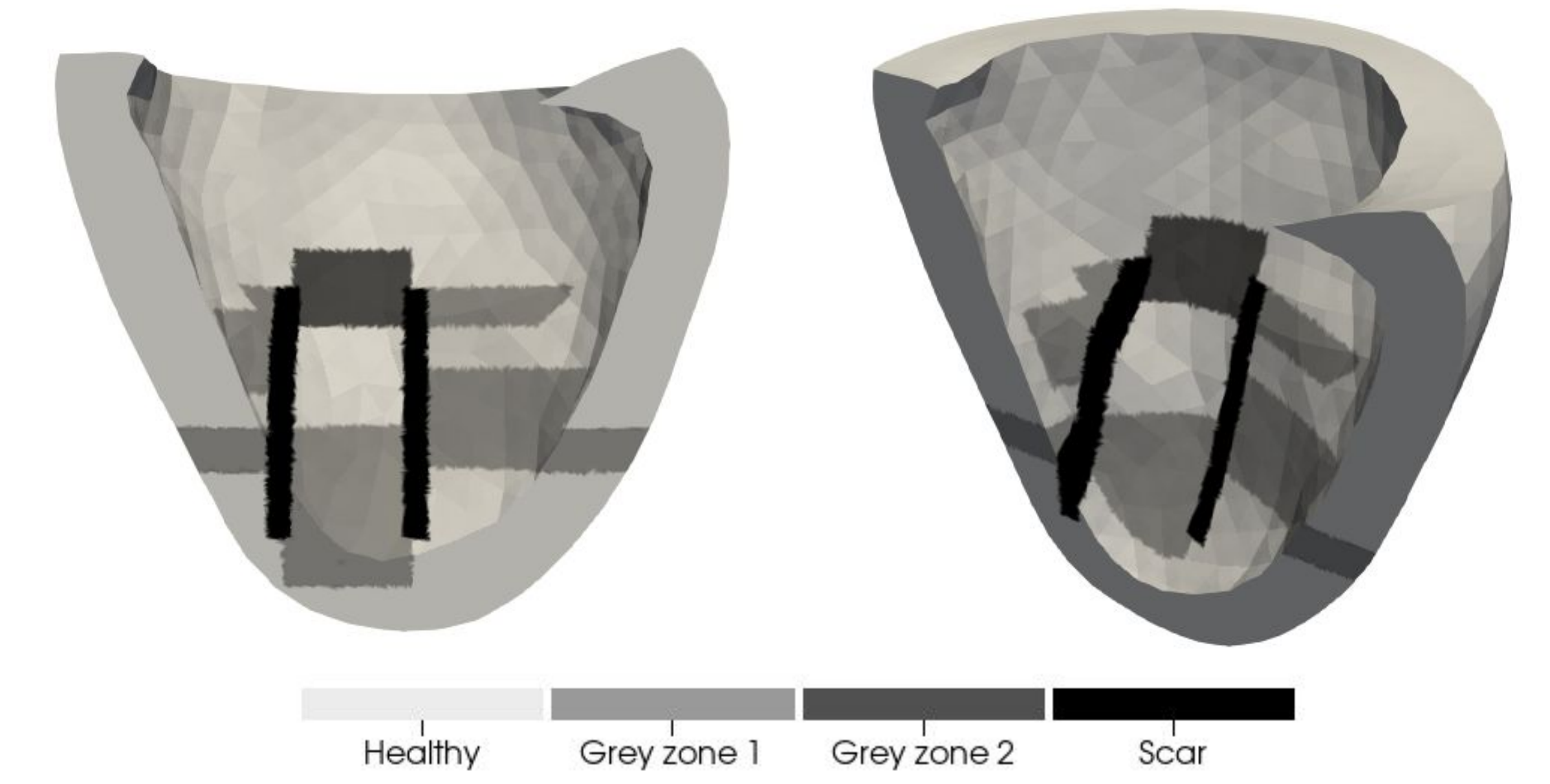} \\
	\EMbaselineVT{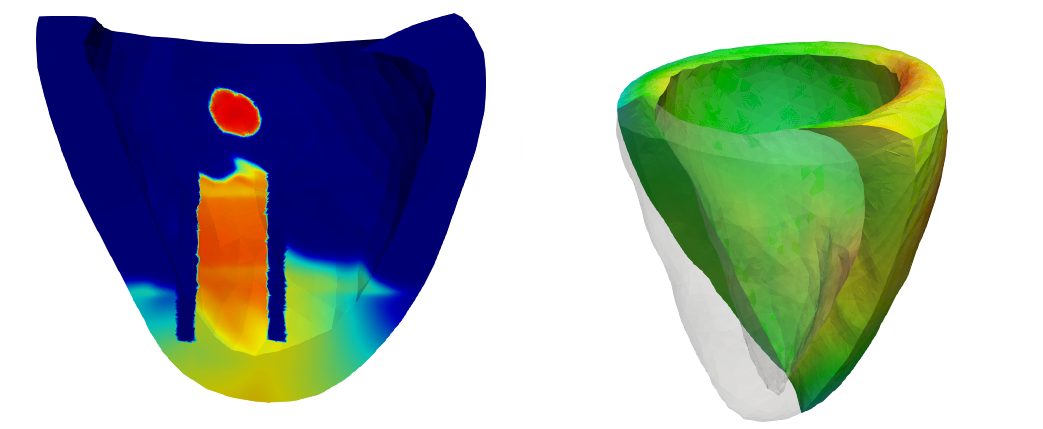}{0.46} \quad
	\EMbaselineVT{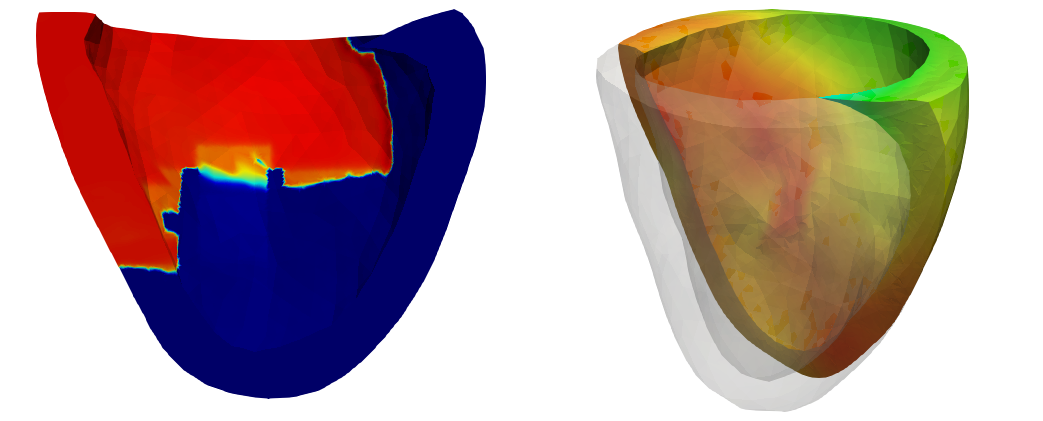}{3.0} \\
	\EMbaselineVT{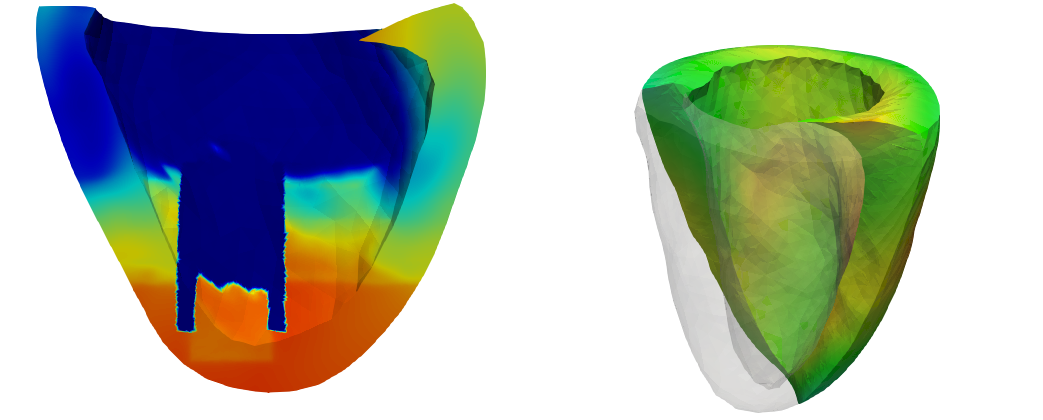}{0.80} \quad
	\EMbaselineVT{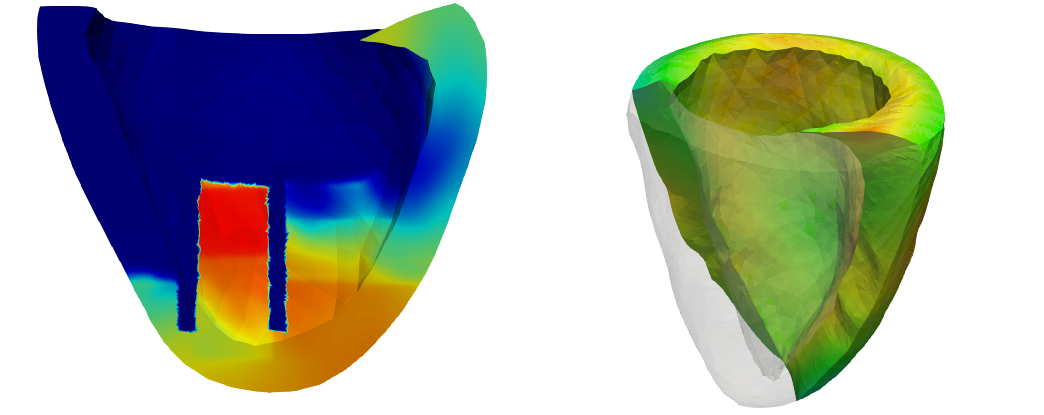}{3.3} \\
	\EMbaselineVT{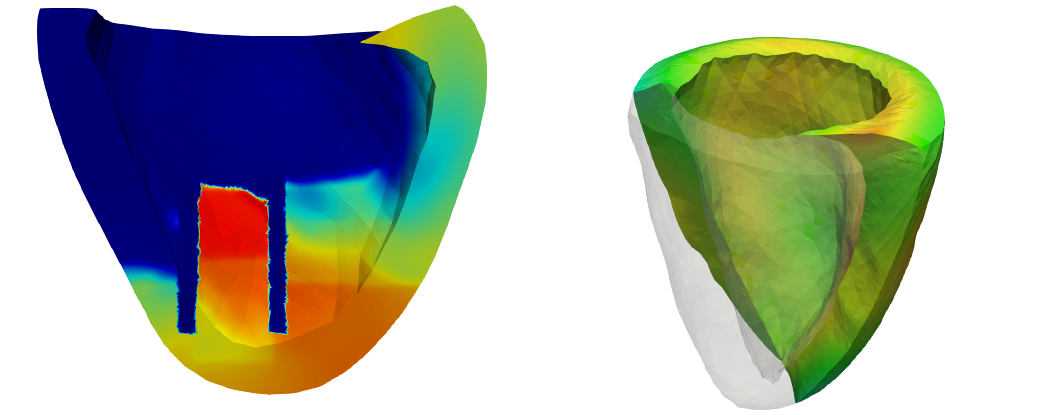}{1.5} \quad
	\EMbaselineVT{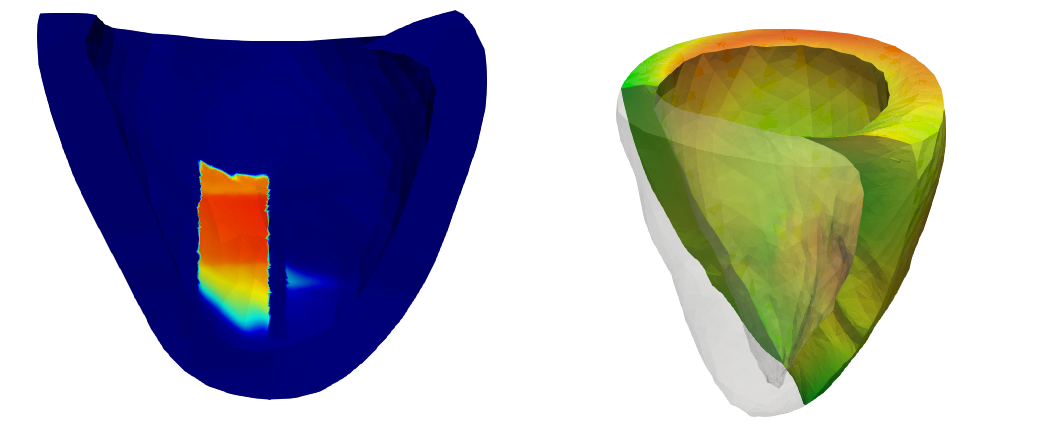}{4.0} \\
  \subfloat[][]{\includegraphics[width=0.5\textwidth]{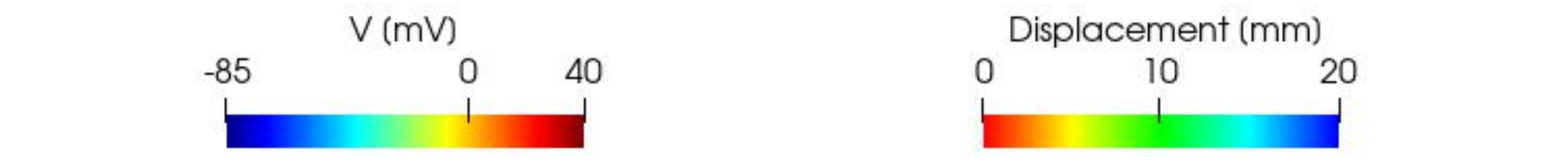}}
  \subfloat[][]{\includegraphics[width=0.5\textwidth]{pictures/EM_VT_none/EM_VT_none_legend.pdf}}

	\caption{Time evolution of the transmembrane potential $V$ (left) and displacement magnitude $|\displ|$ (right) for the Zygote LV with an idealized distribution of ischemic regions.
           Each picture on the right side is warped by the displacement vector $\displ$.
           MEFs are neglected, i.e. we use model $(\Monodomain)$.
					 }
	\label{fig: EM_baseline_VT}
\end{figure}

\ifCIBM \begin{figure}[ht!] \else \begin{figure}[t!] \fi
  \center
  \includegraphics[keepaspectratio, width=0.5\textwidth]{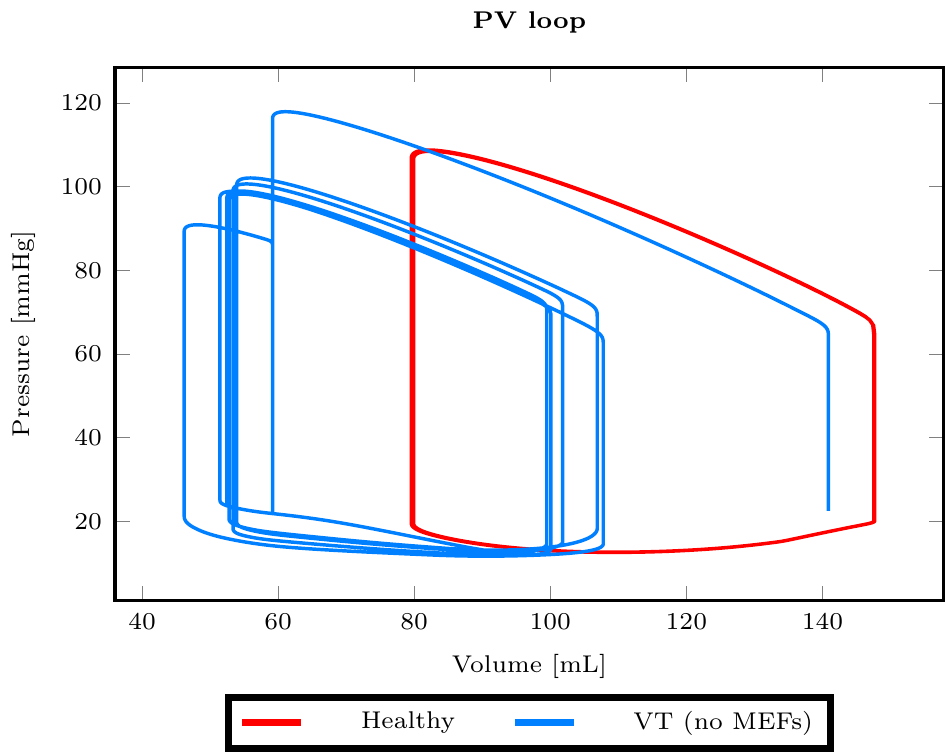}
  \caption{Comparison between a reference healthy PV loop in sinus rhythm (red, Appendix~\ref{app:params}, heartbeat period equal to $\SI{0.8}{\second}$) and the one obtained in the baseline simulation under VT for $t \in [0, 4]$ s (light blue).
           We underline that, differently from \cite{Salvador2021}, here we induce a hemodynamically tolerated VT.}
  \label{fig: EM_baseline_VT_PV_loop}
\end{figure}

We depict in Fig.~\ref{fig: EM_baseline_VT} the evolution of the transmembrane potential and the displacement magnitude for the baseline electromechanical simulation, where all MEFs are fully neglected (model $(\Monodomain)$ in Tab.~\ref{tab: EP_models}).
We induce a sustained VT with a figure-of-eight pattern around the isthmus, which is laterally bordered by scars, that act as conduction blocks.

In Fig.~\ref{fig: EM_baseline_VT_PV_loop}, we compare the Pressure-Volume (PV) loop over different heartbeats for the baseline simulation under VT with a reference healthy PV loop in sinus rhythm obtained by removing scars and grey zones.
We observe that the contractility increases while the stroke volume (SV) decreases.
The ejection fraction (EF) remains approximately the same and we approach a steady state in which the electromechanical function is not impaired.
For these reasons, we conclude that the VT is stable.

VT associated with ischemia are known to disturb the normal isovolumetric processes and to influence the end systolic/diastolic pressure volume relationship (ESPVR/EDPVR).
Simultaneously, a phenomenon called incomplete relaxation may occur, especially when the VT does not leave enough time for the uncoupling of all the actin-myosin bonds between two consecutive contraction phases \cite{Bastos2019}.
The occurrence of this phenomenon is illustrated in Fig.~\ref{fig: incomplete_relaxation}, where we depict the time evolution of the minimum, maximum and average active stress in the computational domain $\Omega_0$ for a reference healthy case in sinus rhythm and the baseline simulation under VT.
Specifically, in the healthy case there is always a time interval between two consecutive heartbeats (precisely, during ventricular diastole) in which the active stress is virtually zero in the LV.
In the VT case, instead, the cardiac muscle is never fully relaxed, thus not allowing the LV to complete its emptying.
All these details are properly captured by our electromechanical model thanks to its biophysical accuracy.

\begin{figure}[t!]
  \center
  \includegraphics[keepaspectratio, width=0.9\textwidth]{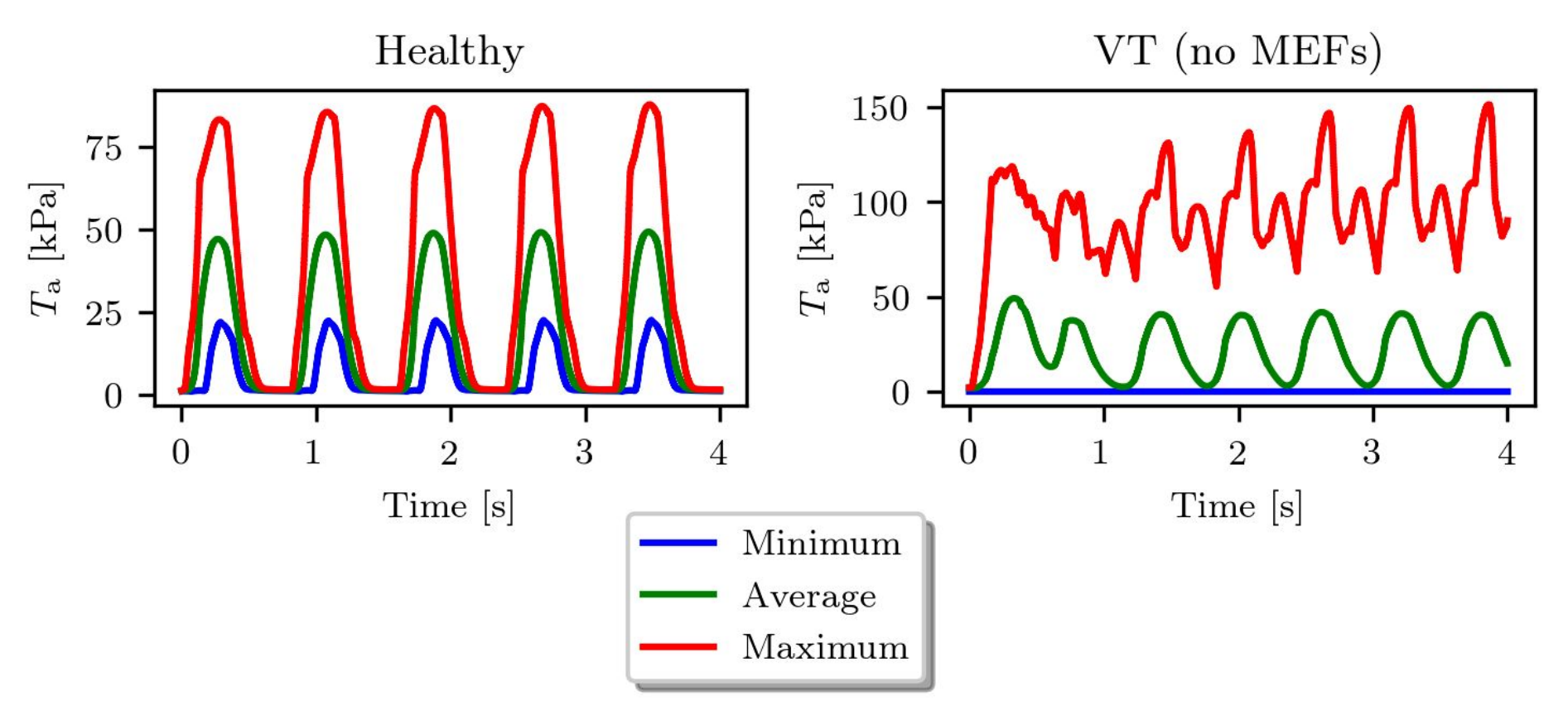}
  \caption{Minimum, average and maximum active tension $\Tens$ over time for a reference healthy case in sinus rhythm (left, Appendix~\ref{app:params}, heartbeat period equal to $\SI{0.8}{\second}$) and the baseline simulation under VT (right).
				   We see that incomplete relaxation occurs during VT.}
  \label{fig: incomplete_relaxation}
\end{figure}

\subsection{Effects of geometry-mediated MEFs}
\label{subsec: gMEFs_effects}

\newcommand{\EMgMEFVT}[2]{
	\subfloat[][$t = \SI{#2}{\second}$]{\includegraphics[width=0.9\textwidth]{pictures/EM_VT_gMEF/EM_VT_gMEF_#1.png}}}

\ifCIBM \begin{figure}[h!] \else \begin{figure}[p] \fi
	\center
	\captionsetup[subfigure]{labelformat=empty}
	\EMgMEFVT{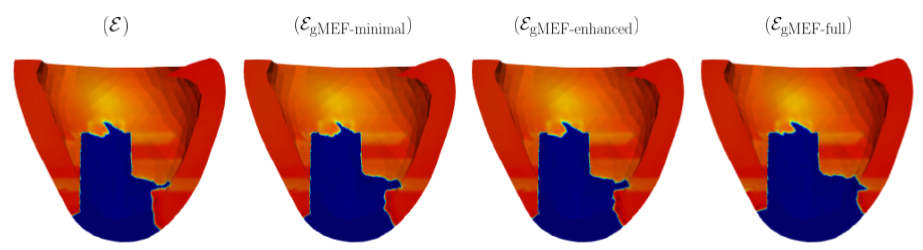}{0.20} \\
	\EMgMEFVT{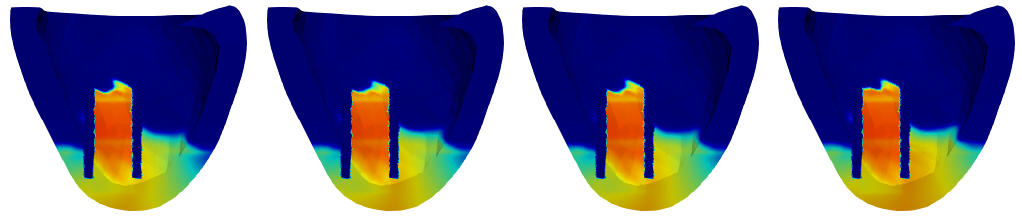}{0.45} \\
	\EMgMEFVT{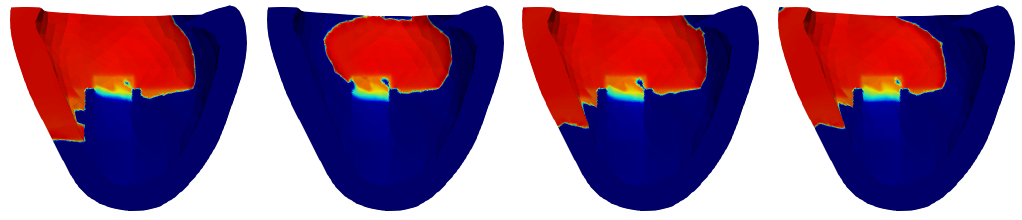}{2.4} \\
	\EMgMEFVT{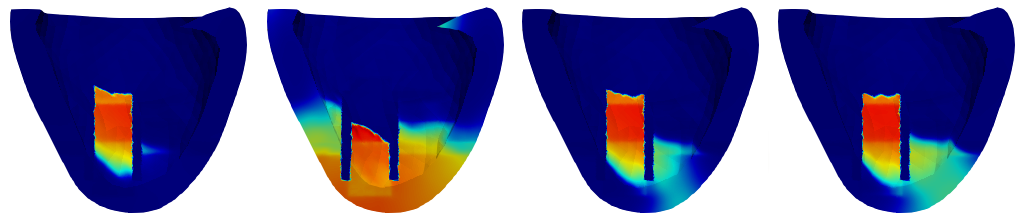}{4.0} \\
  \subfloat[][]{\includegraphics[width=0.6\textwidth]{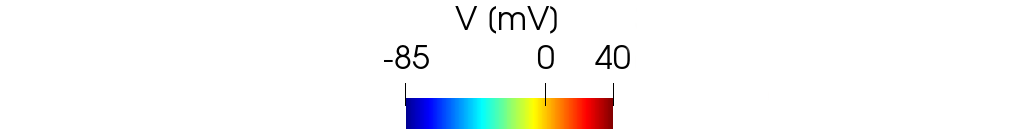}}

	\caption{Comparison among different models for geometry-mediated MEFs in terms of transmembrane potential $V$.
					 }
	\label{fig: EM_gMEF_VT}
\end{figure}

\ifCIBM \begin{figure}[t!] \else \begin{figure} \fi
  \center
  \includegraphics[keepaspectratio, width=0.8\textwidth]{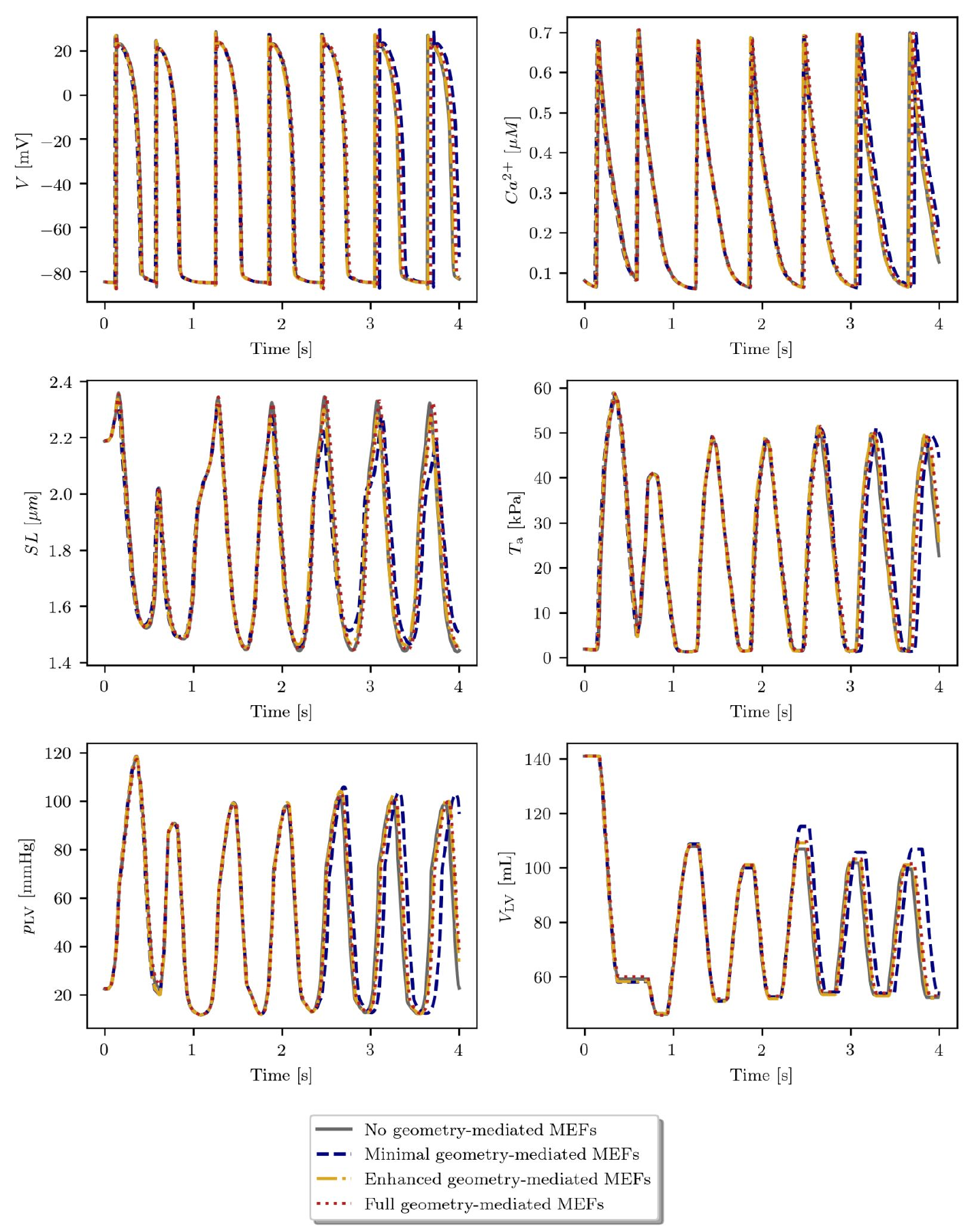}
  \caption{Pointwise values of transmembrane potential $V$, intracellular calcium concentration $\Cai$, sarcomere length $\SL$, active tension $\Tens$, pressure $\PLV$ and volume $\VLV$ over time for $(\Monodomain)$, $(\Monodomain_{\text{gMEF-minimal}})$, $(\Monodomain_{\text{gMEF-enhanced}})$ and $(\Monodomain_{\text{gMEF-full}})$.}
  \label{fig: gMEF_VT}
\end{figure}

\begin{table}
    \begin{center}
        \renewcommand{\arraystretch}{1.5}
        \begin{tabular}{ c | c c c c}
            \toprule
            Model & $(\Monodomain)$ & $(\Monodomain_{\text{gMEF-minimal}})$ & $(\Monodomain_{\text{gMEF-enhanced}})$ & $(\Monodomain_{\text{gMEF-full}})$\\
            \midrule
            BCL & $\SI{0.60}{\second}$ & $\SI{0.65}{\second}$ & $\SI{0.61}{\second}$ & $\SI{0.60}{\second}$\\
            \bottomrule
        \end{tabular}
        \end{center}
        \caption{BCL for different modeling choices in geometry-mediated MEFs.
                 Model $(\Monodomain_{\text{gMEF-minimal}})$ significantly changes BCL with respect to $(\Monodomain)$, $(\Monodomain_{\text{gMEF-enhanced}})$, $(\Monodomain_{\text{gMEF-full}})$.}
        \label{tab:BCL_gMEF}
\end{table}

We consider four different modeling choices for the geometry-mediated MEFs, namely $(\Monodomain)$, $(\Monodomain_{\text{gMEF-minimal}})$, $(\Monodomain_{\text{gMEF-enhanced}})$ and $(\Monodomain_{\text{gMEF-full}})$ in Tab.~\ref{tab: EP_models}, while for the moment we completely neglect the impact of SACs.
Then, we perform four different electromechanical simulations by employing these four different formulations for the monodomain equation.

We illustrate in Fig.~\ref{fig: EM_gMEF_VT} the development of transmembrane potential $V$ over time.
We observe minor differences in action potential propagation among $(\Monodomain)$, $(\Monodomain_{\text{gMEF-enhanced}})$ and $(\Monodomain_{\text{gMEF-full}})$.
These differences, as we can see for $t=\SI{4}{\second}$, are mainly focused on the depolarization wave and occur during VT.
Moreover, by looking at Tab.~\ref{tab:BCL_gMEF}, we notice that the VT BCL is very similar among these three models.
Indeed, the BCL is approximately equal to $\SI{0.60}{\second}$, which is long if compared to more dangerous VT and justifies a stable ventricular excitation.

On the other hand, $(\Monodomain_{\text{gMEF-minimal}})$ entails major changes in VT BCL, which increases from $\SI{0.60}{\second}$ (for model $(\Monodomain)$) to \SI{0.65}{\second}, and conduction velocity, that significantly decreases.

These observations are in agreement with Fig.~\ref{fig: gMEF_VT}, where the electrophysiological, mechanical and hemodynamic variables retrieved in a random point of the computational domain are shifted forward for $(\Monodomain_{\text{gMEF-minimal}})$, while $(\Monodomain)$, $(\Monodomain_{\text{gMEF-enhanced}})$ and $(\Monodomain_{\text{gMEF-full}})$ show a very similar pattern.
This is also motivated by the change in the VT exit site, as we can see from Fig.~\ref{fig: EM_gMEF_VT} for $t = \SI{2.4}{\second}$.
This phenomenon is particularly evident from the plot of sarcomere length over time (Fig.~\ref{fig: gMEF_VT}), which also presents different peak values for $(\Monodomain_{\text{gMEF-minimal}})$ and $t \gtrapprox \SI{2.4}{\second}$.
Finally, we see that wave stability is not affected by geometry-mediated MEFs.
Indeed, the VT always remains hemodynamically stable in the four different cases.

\subsection{Effects of SACs}
\label{subsec: SACs_effects}

\begin{table}
    \begin{center}
        \renewcommand{\arraystretch}{1.5}
        \begin{tabular}{ c | c }
            \toprule
            Model & VT type
            \\
            \midrule
            $(\Monodomain)$ & Stable ($\text{BCL} = \SI{0.60}{\second}$)
            \\
            $(\Monodomain_{\text{SAC}})$, $\Gs = \SI{100}{\per \second}$, $\Vrev = \SI{-70}{\milli \volt}$ & Stable ($\text{BCL} = \SI{0.60}{\second}$)
            \\
            $(\Monodomain_{\text{SAC}})$, $\Gs = \SI{100}{\per \second}$, $\Vrev = \SI{-35}{\milli \volt}$ & Stable ($\text{BCL} = \SI{0.60}{\second}$)
            \\
            $(\Monodomain_{\text{SAC}})$, $\Gs = \SI{100}{\per \second}$, $\Vrev = \SI{0}{\milli \volt}$ & Unstable ($\text{BCL}_\text{avg} = \SI{0.50}{\second}$)
            \\
            $(\Monodomain_{\text{SAC}})$, $\Gs = \SI{50}{\per \second}$, $\Vrev = \SI{0}{\milli \volt}$ & Stable ($\text{BCL} = \SI{0.60}{\second}$)
            \\
            \bottomrule
        \end{tabular}
        \end{center}
        \caption{VT classification for different SACs parametrizations.
                 The unstable VT has a BCL that ranges from $\SI{0.43}{\second}$ to $\SI{0.58}{\second}$.}
        \label{tab:BCL_SAC}
\end{table}

In this section, we fully neglect the effects of geometry-mediated MEFs and we focus on different parametrizations for SACs in terms of $\Gs$ and $\Vrev$.
We also compare the outcomes of $(\Monodomain_{\text{SAC}})$ with the ones of $(\Monodomain)$ to outline similarities and differences.

We notice from Fig.~\ref{fig: EM_SAC_VT} that both APD and wave stability are affected by SACs parametrizations.
Indeed, by combining the 3D information with the pointwise evaluations of Fig.~\ref{fig: SAC_100_VT}, we discover that there is one choice of the parameters, namely $\Gs = \SI{100}{\per \second}$ and $\Vrev = \SI{0}{\milli \volt}$, that converts the VT from stable to unstable.
The instability derives from an extra stimulus that is completely driven by contraction, which occurs in the superior-right part of the ventricle.
This extra stimulus changes the VT morphology, along with its BCL, which is not the same over time.

By considering data of Tab.~\ref{tab:BCL_SAC}, we observe that the VT BCL remains the same when there is no stability transition.
Moreover, from the hemodynamic perspective, the onset of different types of arrhythmias is driven by the combined effects of $\Gs$ and $\Vrev$.
Indeed, in Fig.~\ref{fig: SAC_100_VT} we see that given $\Gs$, different $\Vrev$ may change wave stability.
On the other hand, in Fig.~\ref{fig: EM_SAC_50_100_VT_pressure_volume} we show that different $\Gs$ may affect wave stability, with $\Vrev$ fixed a priori.

\newcommand{\EMSACVT}[2]{
	\subfloat[][$t = \SI{#2}{\second}$]{\includegraphics[width=0.9\textwidth]{pictures/EM_VT_SAC_100/EM_VT_SAC_100_#1.png}}}

\ifCIBM \begin{figure}[h!] \else \begin{figure}[p] \fi
	\center
	\captionsetup[subfigure]{labelformat=empty}
	\EMSACVT{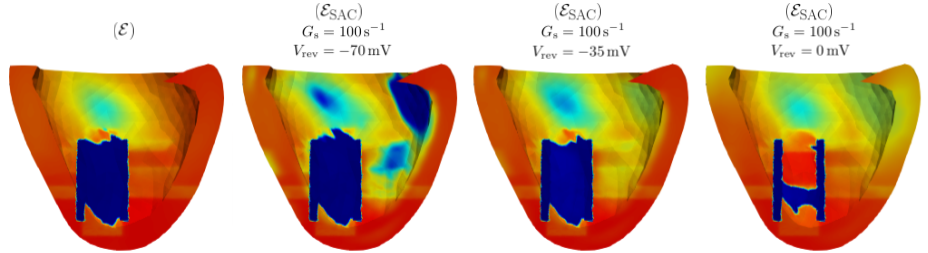}{0.25} \\
	\EMSACVT{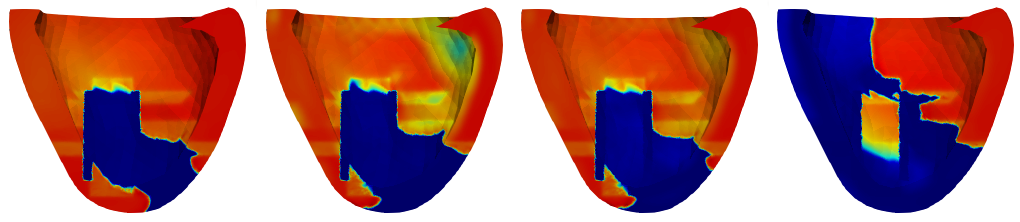}{2.5} \\
	\EMSACVT{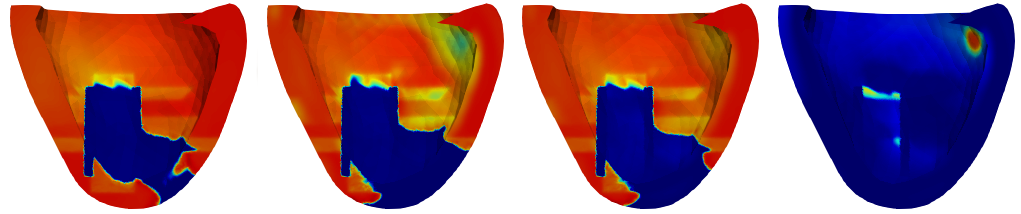}{3.1} \\
	\EMSACVT{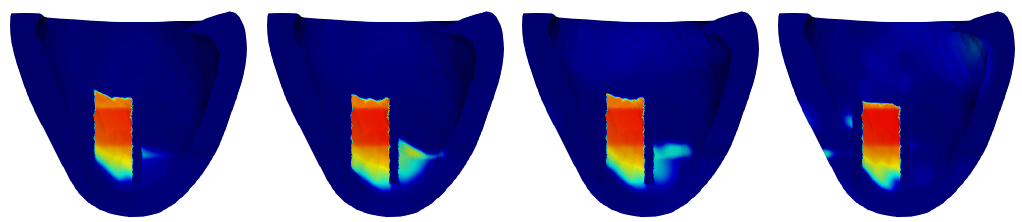}{4.0} \\
  \subfloat[][]{\includegraphics[width=0.6\textwidth]{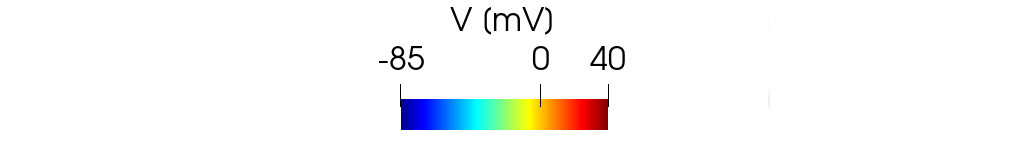}}

	\caption{Comparison between model $(\Monodomain)$ and model $(\Monodomain_{\text{SAC}})$, for different values of $\Vrev$ ($\Gs = \SI{100}{\per \second}$).
					 }
	\label{fig: EM_SAC_VT}
\end{figure}

\ifCIBM \begin{figure}[h!] \else \begin{figure} \fi
  \center
  \includegraphics[keepaspectratio, width=0.8\textwidth]{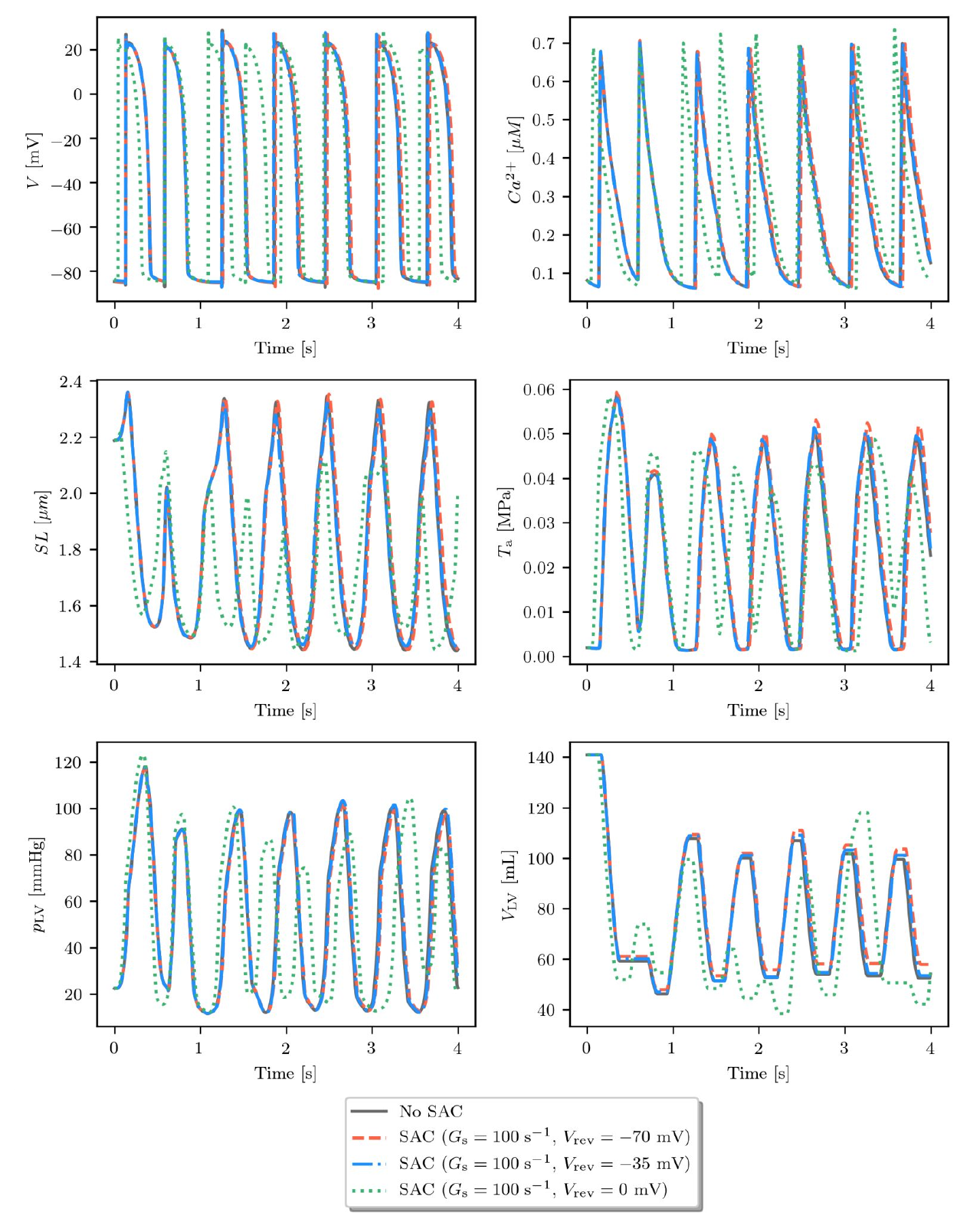}
  \caption{Pointwise values of transmembrane potential $V$, intracellular calcium concentration $\Cai$, sarcomere length $\SL$, active tension $\Tens$, pressure $\PLV$ and volume $\VLV$ over time for $(\Monodomain)$ and $(\Monodomain_{\text{SAC}})$ with different choices of $\Vrev$ ($\Gs = \SI{100}{\per \second}$).}
  \label{fig: SAC_100_VT}
\end{figure}

\ifCIBM \begin{figure}[h!] \else \begin{figure} \fi
  \center
  \includegraphics[keepaspectratio, width=0.8\textwidth]{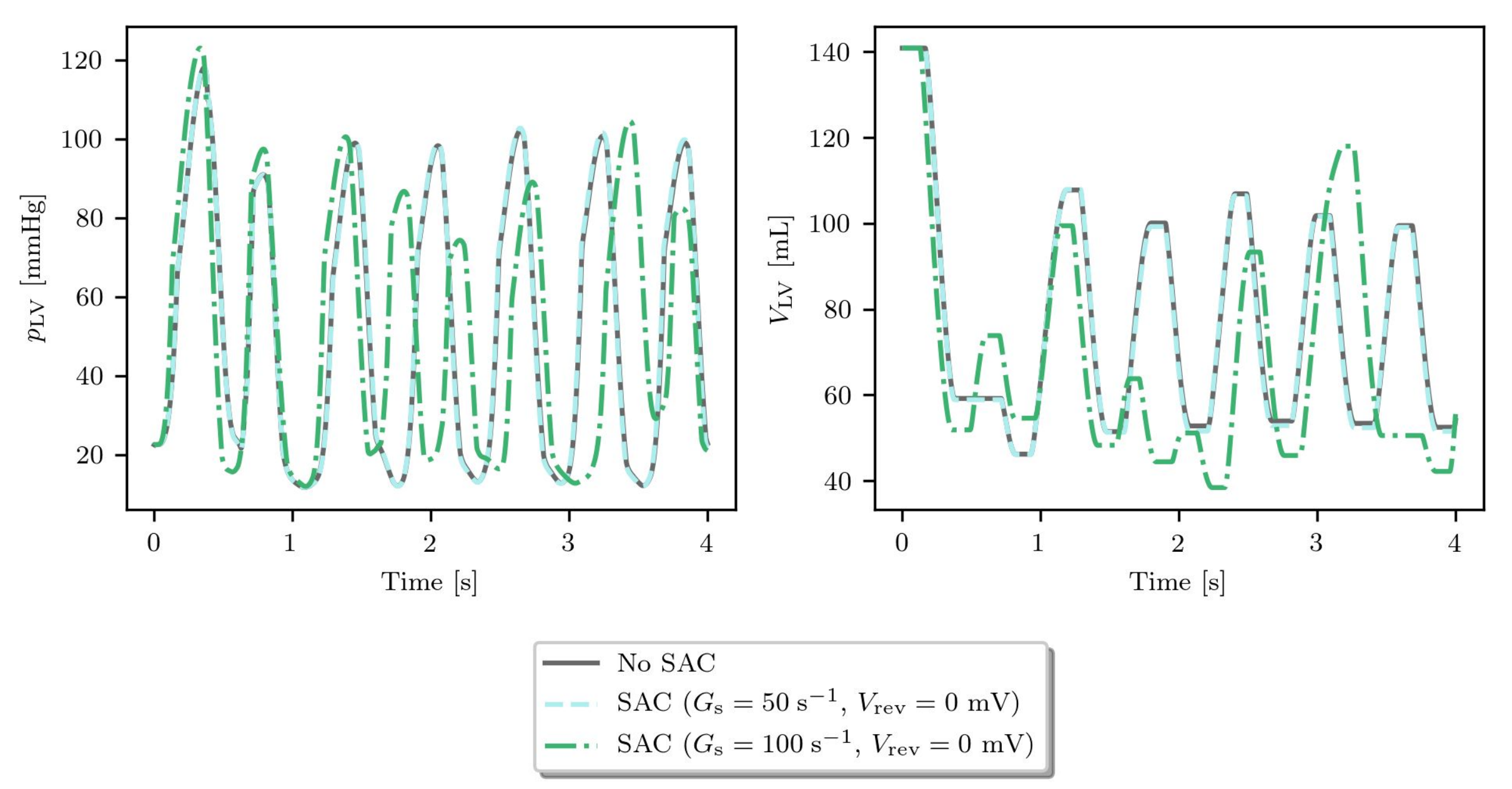}
  \caption{Pointwise values of pressure $\PLV$ and volume $\VLV$ over time for $(\Monodomain)$ and $(\Monodomain_{\text{SAC}})$ with different choices of $\Gs$ ($\Vrev = \SI{0}{\milli \volt}$)}
  \label{fig: EM_SAC_50_100_VT_pressure_volume}
\end{figure}

\subsection{Combined effects of geometry-mediated MEFs and SACs}
\label{subsec: MEFs_effects}

\newcommand{\EMgMEFSACVT}[2]{
	\subfloat[][$t = \SI{#2}{\second}$]{\includegraphics[width=0.75\textwidth]{pictures/EM_VT_gMEF_SAC_100/EM_VT_gMEF_SAC_100_#1.png}}}

\ifCIBM \begin{figure}[h!] \else \begin{figure}[p] \fi
	\center
	\captionsetup[subfigure]{labelformat=empty}
	\EMgMEFSACVT{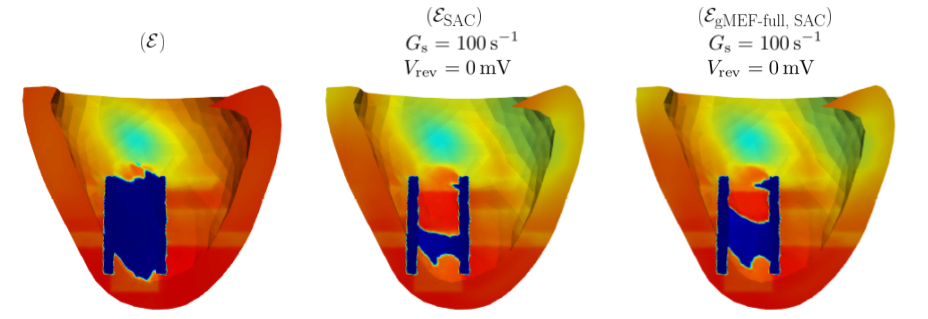}{0.25} \\
	\EMgMEFSACVT{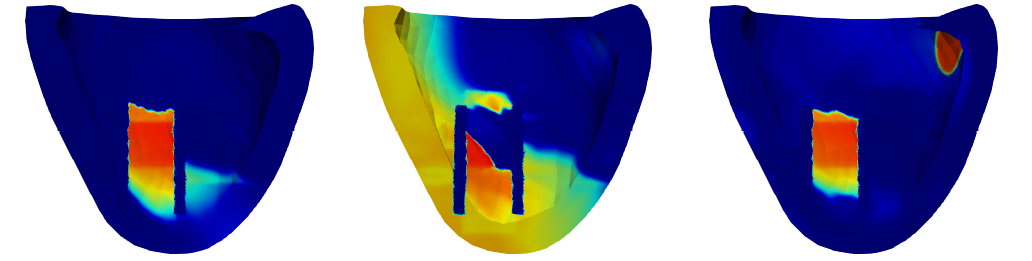}{2.79} \\
	\EMgMEFSACVT{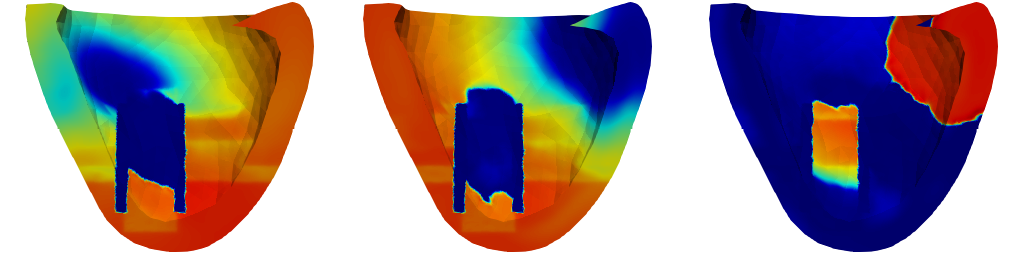}{3.80} \\
	\EMgMEFSACVT{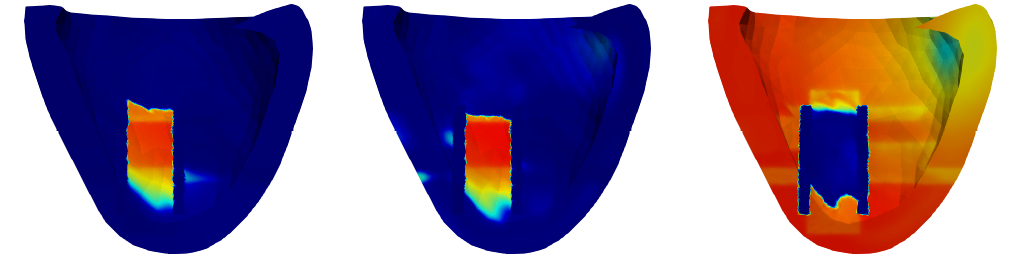}{4.0} \\
  \subfloat[][]{\includegraphics[width=0.6\textwidth]{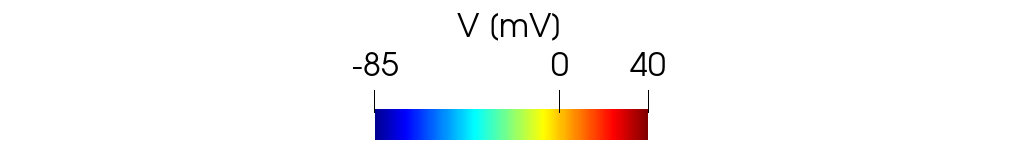}}

	\caption{Comparison among models $(\Monodomain)$, $(\Monodomain_{\text{SAC}})$ and $(\Monodomain_{\text{gMEF-full, SAC}})$.
					 }
	\label{fig: EM_gMEF_SAC_VT}
\end{figure}

\begin{figure}[t!]
  \center
  \includegraphics[keepaspectratio, width=0.8\textwidth]{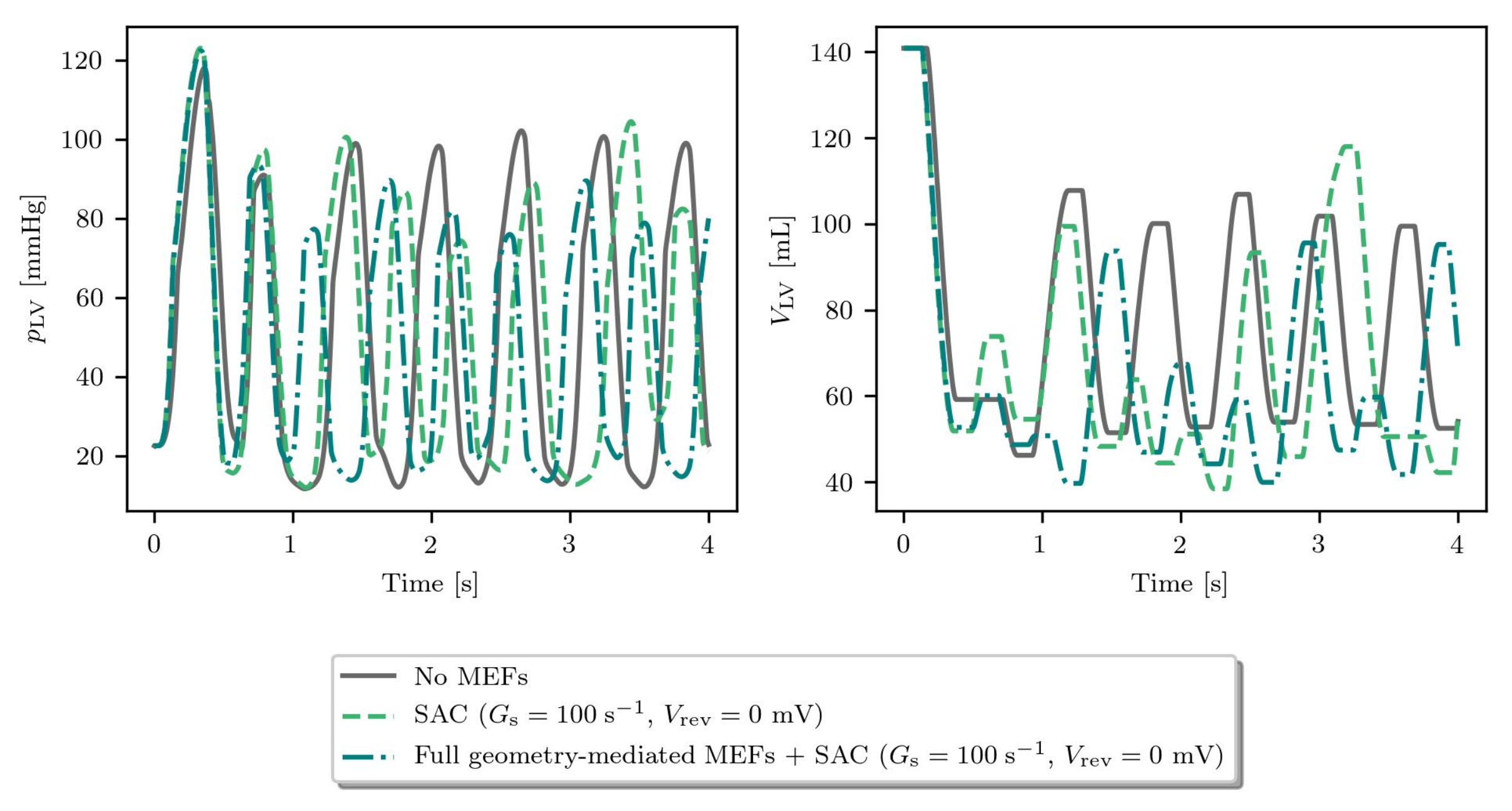}
  \caption{Pointwise values of pressure $\PLV$ and volume $\VLV$ over time for $(\Monodomain)$, $(\Monodomain_{\text{SAC}})$ ($\Gs = \SI{100}{\per \second}$, $\Vrev = \SI{0}{\milli \volt}$) and $(\Monodomain_{\text{gMEF-full, SAC}})$ ($\Gs = \SI{100}{\per \second}$, $\Vrev = \SI{0}{\milli \volt}$).}
  \label{fig: EM_gMEF_SAC_100_VT_pressure_volume}
\end{figure}

\ifCIBM \begin{figure}[h!] \else \begin{figure}[p!] \fi
  \center
  \includegraphics[keepaspectratio, width=0.9\textwidth]{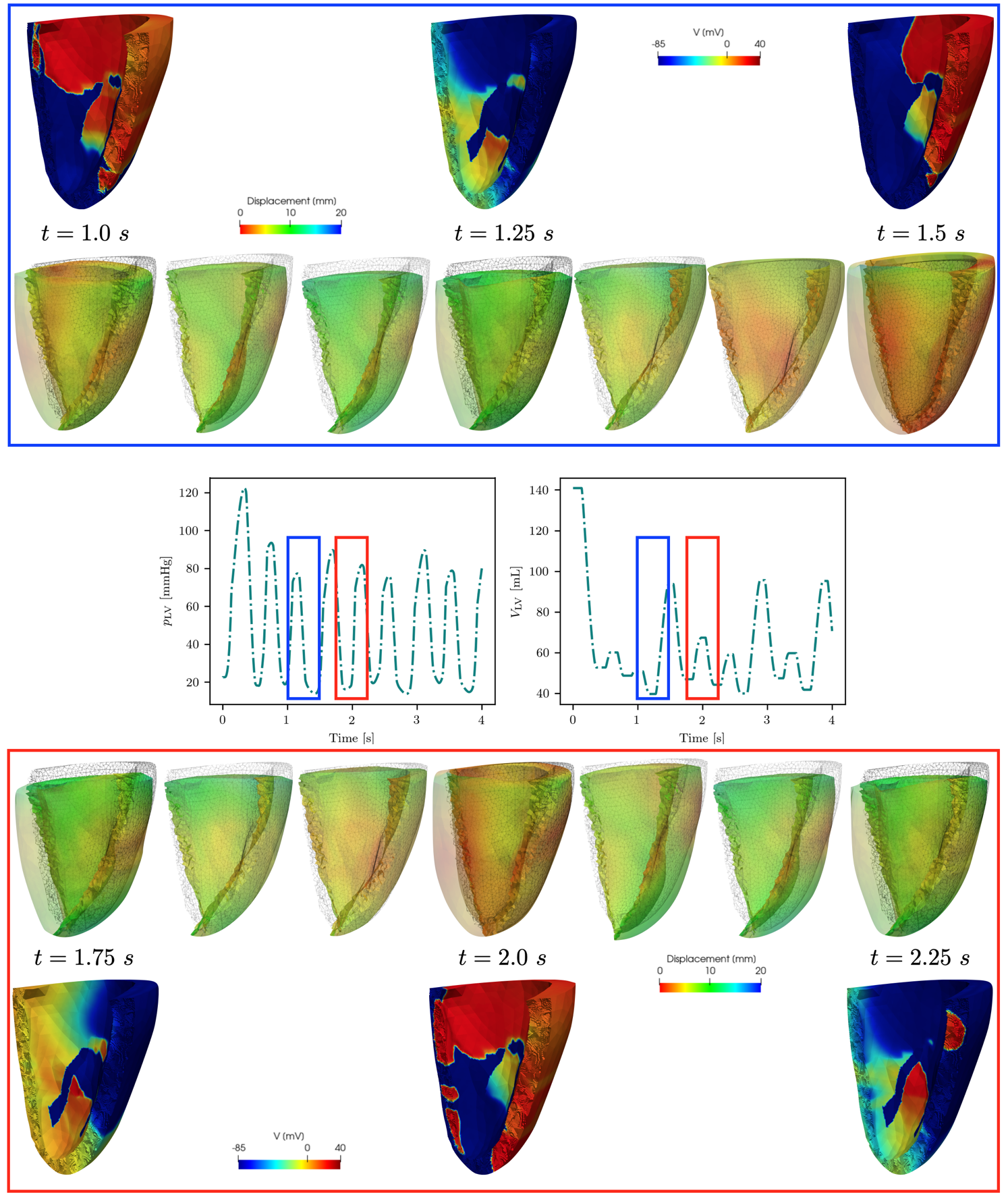}
  \caption{Coupled effects of electrophysiology, mechanics and hemodynamics for the numerical simulation with model $(\Monodomain_{\text{gMEF-full, SAC}})$.
	         The extra stimuli in the upper right part of the LV, which is driven by SACs, activate the LV electrophysiologically and mechanically.
					 This has a direct impact on both pressure and volume transients, which in turn have an effect on the electromechanical behavior of the LV.}
  \label{fig: gMEF_full_SAC}
\end{figure}

\begin{table}
    \begin{center}
        \renewcommand{\arraystretch}{1.5}
        \begin{tabular}{ c | c }
            \toprule
            Model & VT type
            \\
            \midrule
            $(\Monodomain)$ & Stable ($\SI{0.60}{\second}$)
            \\[1mm]
            $(\Monodomain_{\text{SAC}})$, $\Gs = \SI{100}{\per \second}$, $\Vrev = \SI{0}{\milli \volt}$ & Unstable ($\text{BCL}_\text{avg} = \SI{0.50}{\second}$)
            \\
            $(\Monodomain_{\text{gMEF-full, SAC}})$, $\Gs = \SI{100}{\per \second}$, $\Vrev = \SI{0}{\milli \volt}$ & Unstable ($\text{BCL}_\text{avg} = \SI{0.47}{\second}$)
            \\
            \bottomrule
        \end{tabular}
        \end{center}
        \caption{VT classification for combinations of geometry-mediated MEFs and nonselective SACs.
                 The unstable VT related to $(\Monodomain_{\text{SAC}})$ has a BCL that ranges from $\SI{0.43}{\second}$ to $\SI{0.58}{\second}$.
                 The unstable VT related to $(\Monodomain_{\text{gMEF-full, SAC}})$ has a BCL that ranges from $\SI{0.44}{\second}$ to $\SI{0.50}{\second}$.}
        \label{tab:BCL_gMEF_SAC_100}
\end{table}

We briefly evaluate the combined effects of geometry-mediated MEFs and nonselective SACs.
Once SACs parametrization is fixed, we notice that switching between no formulation (i.e. $(\Monodomain_{\text{SAC}})$) to full formulation (i.e. $(\Monodomain_{\text{gMEF-full, SAC}})$) of geometry-mediated MEFs entails significant differences.
In particular, from Fig.~\ref{fig: EM_gMEF_SAC_VT} we see that $(\Monodomain_{\text{gMEF-full, SAC}})$ triggers the extra stimuli faster than $(\Monodomain_{\text{SAC}})$.
As we can notice from Fig.~\ref{fig: EM_gMEF_SAC_100_VT_pressure_volume}, the VT remains unstable but both pressure and volume traces over time are very different from each other for $(\Monodomain_{\text{SAC}})$ and $(\Monodomain_{\text{gMEF-full, SAC}})$.
From Tab.~\ref{tab:BCL_gMEF_SAC_100}, we infer that the VT BCL for $(\Monodomain_{\text{gMEF-full, SAC}})$ is lower than the one of $(\Monodomain_{\text{SAC}})$.
This potentially defines a more dangerous VT.
Finally, in Fig.~\ref{fig: gMEF_full_SAC} we highlight the joint contributions of electrophysiology, mechanics and hemodynamics in the $(\Monodomain_{\text{gMEF-full, SAC}})$ coupled model.

	\section{Discussion}
\label{sec: discussion}

We investigate how several types of geometry-mediated MEFs and different parametrizations for nonselective SACs affect the electric and hemodynamic stability of VT.
In particular, we focus on the sustainment and the morphology of VT macro-reentrant circuits and blood supply, which is analyzed by means of PV loops.

Differently from previous studies \cite{ColliFranzone2017, Timmermann2017AnIA, Panfilov2010}, we keep into account tissue heterogeneity by introducing an idealized distribution of scars, grey zones and non-remodeled regions over the myocardium.
We also consider a more sophisticated coupled mathematical model that embraces electrophysiology, activation, mechanics and cardiovascular fluid dynamics to analyze these mechano-electric couplings.
To the best of our knowledge, this framework enables to study for the first time the hemodynamic effects of both geometry-mediated MEFs and nonselective SACs on VT by means of electromechanical simulations, showing strengths and weaknesses of some simplifying modeling choices that are commonly made when simulating this type of phenomenon \cite{Chapelle2009, Collin2019, levrero2020sensitivity, quarteroni2017integrated}.
The coupling between the 3D electromechanical model and the 0D circulation model permits to identify the hemodynamic nature of the VT.
With our approach, we discriminate between stable and unstable VT, which might result in hemodynamically tolerated or not tolerated VT \cite{Salvador2021}.

We see that if a VT is triggered by a certain stimulation protocol and by neglecting all MEFs, the very same pacing protocol induces a VT for all possible combinations of MEFs \cite{Salvador2021}.

In our numerical simulations, we do not observe a significant impact of geometry-mediated MEFs on the induction and sustainment of the VT. 
Most of the modeling choices for geometry-mediated MEFs, namely $(\Monodomain)$, $(\Monodomain_{\text{gMEF-enhanced}})$ and $(\Monodomain_{\text{gMEF-full}})$, present very similar VT BCLs and conduction velocities, while showing a few differences in the depolarization wave \cite{ColliFranzone2017}.
On the other hand, $(\Monodomain_{\text{gMEF-minimal}})$ manifests major differences in the VT BCL, which increases, and its exit site with respect to $(\Monodomain)$, $(\Monodomain_{\text{gMEF-enhanced}})$ and $(\Monodomain_{\text{gMEF-full}})$, as shown in \cite{Salvador2021} for a patient-specific unstable VT.
This can be justified by the simplifications introduced in $(\Monodomain_{\text{gMEF-minimal}})$, while moving towards $(\Monodomain_{\text{gMEF-full}})$ we almost totally recover the behavior observed in $(\Monodomain)$.
Therefore, the minimal MEFs modeling choice might lead to biased results which are in favor of less severe VT.



We observe that nonselective SACs may affect the hemodynamic nature of the VT, as they may induce EADs or DADs, which lead to ectopic foci that reactivate the LV \cite{Hazim2021}.
These extra stimuli are generally located in the regions of the myocardium in which there is a transition between scar and border zone or between border zone and non-remodeled areas, where high stretches are likely to be present \cite{Jie2010}.
According to the specific combination of $\Gs$ and $\Vrev$, nonselective SACs affects both APD and AP resting values \cite{Kohl2004}.
We remark that such spontaneous arrhythmias triggered by myocardial stretches cannot be assessed in electrophysiological simulations, where the mechanical behavior is neglected.



Finally, in our study we also investigate the combined effects of geometry-mediated MEFs and nonselective SACs.
Significant differences between $(\Monodomain_{\text{gMEF, full, SAC}})$ and $(\Monodomain_{\text{SAC}})$ are observed when geometry-mediated MEFs are combined with a parametrization of SACs that entails extra stimuli.
Specifically, the VT BCL with $(\Monodomain_{\text{gMEF, full, SAC}})$ is lower than the one of $(\Monodomain_{\text{SAC}})$ because the extra stimuli driven by SACs is triggered more often in the former case.
This completely changes the pressure-volume dynamics of the VT, whose stability is however still not affected by the geometry-mediated MEFs.

	\section{Limitations}
\label{sec: limitations}

There are other types of MEFs that could be investigated in future works.
Among them, an important role is certainly played by the mechanical effects mediated by fibroblasts in the extracellular matrix, $\Cai$ buffers handling, alterations in transmembrane capacitance $\EPCm$ due to local stretch and ions selective SACs \cite{Kohl2013}.
Indeed, modeling different cellular processes in cardiac electrophysiology might be of interest to further shed light on the underlying mechanisms of arrhythmias.
Nevertheless, we have provided a broad study of geometry-mediated MEFs and nonselective SACs by combining electrophysiological, mechanical and hemodynamic observations.

Furthermore, MEFs should be properly quantified in patient-specific cases.
Currently, our model for nonselective SACs permits to put one single value of both $\Gs$ and $\urev$ for the whole geometry, which is unlikely to happen in realistic scenarios.
A precise space assessment of SACs combined with our mathematical framework could have significant clinical implications.


	\section{Conclusions}
\label{sec: conclusions}

We studied the effects of geometry-mediated MEFs and nonselective SACs on a realistic LV geometry endowed with an idealized distribution of infarct and peri-infarct zones.
We performed numerical simulations of cardiac electromechanics coupled with closed-loop cardiovascular circulation under VT.

Our electromechanical framework allows for the hemodynamic classification of the VT, which can be either stable or unstable, and permits to capture mechanically relevant indications under VT, such as the incomplete relaxation of sarcomeres.
Furthermore, by combining electrophysiology, activation, mechanics and hemodynamics, we observed several differences on the morphology of the VT with respect to electrophysiological simulations.
In particular, geometry-mediated MEFs do not affect VT stability but may alter the VT BCL, along with its exit site \cite{Panfilov2010}.
On the other hand, the recruitment of SACs may generate extra stimuli, which may change VT stability.
These extra stimuli are driven by myocardial contraction and are induced by changes in the APD or in the resting value of the transmembrane potential \cite{Jie2010}.
We conclude that both geometry-mediated MEFs and nonselective SACs define important contributions in electromechanical models with hemodynamic coupling, especially when numerical simulations under arrhythmia are carried out.

	\section*{Acknowledgements}
MS, FR, SP, LD and AQ acknowledge the European Research Council (ERC) under the European Union's Horizon 2020 research and innovation programme (grant agreement No 740132, iHEART - An Integrated Heart Model for the simulation of the cardiac function, P.I. Prof. A. Quarteroni).
NT acknowledges the National Institutes of Health (grants R01HL142496 and R01HL142893) and a Leducq Foundation grant.
\begin{center}
	\raisebox{-.5\height}{\includegraphics[width=.15\textwidth]{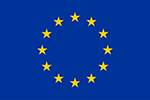}}
	\hspace{2cm}
	\raisebox{-.5\height}{\includegraphics[width=.15\textwidth]{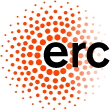}}
\end{center}

	\begin{appendices}
		\section{Model parameters} \label{app:params}

We provide the full list of parameters adopted for the numerical simulations of this paper.
Specifically, Tab.~\ref{tab:parameters electrophysiology} contains the parameters related to the electrophysiological model, Tab.~\ref{tab:parameters sarcomere} those related to the sarcomere model RDQ20-MF, Tab.~\ref{tab:parameters mechanics} the parameters of the mechanical model and, finally, Tab.~\ref{tab:parameters circulation} contains the parameters associated with the circulation model.
For the TTP06 model, we adopt the parameters reported in the original paper (for endocardial cells) \cite{ten2006alternans}, with the only difference that we rescale the intracellular calcium concentration $\Cai$ by a factor of $\CaiFactor$ to get more physiological values \cite{coppini2013late}.

\begin{table}[h!]
	\centering
	\begin{tabular}{lrr|lrr}
		\toprule
		Variable & Value & Unit & Variable & Value & Unit \\
		\midrule
		\multicolumn{3}{l|}{\textbf{Conductivity tensor}} & \multicolumn{3}{l}{\textbf{Applied current}} \\
		$\sigma_{\text{l}}$ & \num{0.7643e-4} & \si{\meter\squared\per\second}   &
		$\EPIappReducedMax$  & \num{17} & \si{\volt\per\second}   \\
		$\sigma_{\text{t}}$ & \num{0.3494e-4} & \si{\meter\squared\per\second}   &
		$\EPIappDuration$    & \num{3e-3} & \si{\second}   \\
		$\sigma_{\text{n}}$ & \num{0.1125e-4} & \si{\meter\squared\per\second}   &
		\multicolumn{3}{l}{\textbf{Calcium rescaling}} \\
		& & &
		$\CaiFactor$ & \num{0.48} & - \\
		\bottomrule
	\end{tabular}
	\caption{Parameters of the electrophysiological model.}
	\label{tab:parameters electrophysiology}
\end{table}

\begin{table}[h!]
	\centering
	\begin{tabular}{lrr|lrr}
		\toprule
		Variable & Value & Unit & Variable & Value & Unit \\
		\midrule
		\multicolumn{3}{l|}{\textbf{Regulatory units steady-state}}              & \multicolumn{3}{l}{\textbf{Crossbridge cycling}}         \\
			$\mu$          & \num{10}       & -                                   &    $\mu_{\fP}^{0}$& \num{32.255}   & $\si{\per\second}$  \\
			$\gamma$       & \num{30}       & -                                   &    $\mu_{\fP}^{1}$& \num{0.768}    & $\si{\per\second}$  \\
			$Q$            & \num{2}        & -                                   &    $r_0$          & \num{134.31}   & $\si{\per\second}$  \\
			$\KdZero$      & \num{0.4}      & $\si{\micro\molar}$                 &    $\alpha$       & \num{25.184}   & -                   \\
			$\KdAlpha$     & \num{-0.2083}  & $\si{\micro\molar\per\micro\meter}$ & \multicolumn{3}{l}{\textbf{Upscaling}}                   \\
		\multicolumn{3}{l|}{\textbf{Regulatory units kinetics}}                  &    $\aXB$         & \num{160}      & $\si{\mega\pascal}$ \\
			$\Koff$        & \num{40}       & $\si{\per\second}$                  &    $\SL_0$        & \num{1.9}      & $\si{\micro\meter}$ \\
			$\Kbasic$      & \num{8}        & $\si{\per\second}$                  &                   &                &                     \\
		\bottomrule
	\end{tabular}
	\caption{Parameters of the sarcomere model RDQ20-MF (for the definition of the parameters, see \cite{regazzoni2020biophysically}).}
	\label{tab:parameters sarcomere}
\end{table}

\begin{table}[h!]
	\centering
	\begin{tabular}{lrr|lrr}
		\toprule
		Variable & Value & Unit & Variable & Value & Unit \\
		\midrule
		\multicolumn{3}{l|}{\textbf{Constitutive law}} & \multicolumn{3}{l}{\textbf{Boundary conditions}} \\
		$B$   & \num{50e3}           & \si{\pascal}   &
		$\BCmecKepiN$      & \num{2e5}       & \si{\pascal\per\meter}   \\
		$C$    & \num{0.88e3}        & \si{\pascal}   &
		$\BCmecKepiT$      & \num{2e4}       & \si{\pascal\per\meter}   \\
		$b_{\text{ff}}$    & 8          & $-$   &
		$\BCmecCepiN$      & \num{2e4}       & \si{\pascal\second\per\meter}   \\
		$b_{\text{ss}}$    & 6        & $-$   &
		$\BCmecCepiT$   &  \num{2e3}      & \si{\pascal\second\per\meter}   \\
		$b_{\text{nn}}$    & 3          & $-$   &
		\multicolumn{3}{l}{\textbf{Tissue density}} \\
		$b_{\text{fs}}$    & 12         & $-$   &
		$\rho_\text{s}$    & $10^3$         & \si{\kilogram \per \cubic\meter} \\
		$b_{\text{fn}}$    & 3         & $-$   &
		\\
		$b_{\text{sn}}$    & 3         & $-$   &
		\\
		\bottomrule
	\end{tabular}
	\caption{Parameters of the mechanical model.}
	\label{tab:parameters mechanics}
\end{table}

\begin{table}[h!]
	\centering
	\begin{tabular}{lrr|lrr}
		\toprule
		Variable & Value & Unit & Variable & Value & Unit \\
		\midrule
		\multicolumn{3}{l|}{\textbf{External circulation}} & \multicolumn{3}{l}{\textbf{Cardiac chambers}} \\
		$\RarSYS$   & 0.64           & \si{\mmHg \second \per \milli\liter}   &
		$\EpLA$      & 0.18       & \si{\mmHg \per \milli\liter}   \\
		$\RarPUL$    & 0.032116        & \si{\mmHg \second \per \milli\liter}   &
		$\EpRA$      & 0.07       & \si{\mmHg \per \milli\liter}   \\
		$\RvnSYS$    & 0.32          & \si{\mmHg \second \per \milli\liter}   &
		$\EpRV$      & 0.05       & \si{\mmHg \per \milli\liter}   \\
		$\RvnPUL$    & 0.035684        & \si{\mmHg \second \per \milli\liter}   &
		$\EaMaxLA$   & 0.07       & \si{\mmHg \per \milli\liter}   \\
		$\CarSYS$    & 1.2          & \si{\milli\liter \per \mmHg}   &
		$\EaMaxRA$   & 0.06       & \si{\mmHg \per \milli\liter}   \\
		$\CarPUL$    & 10.0         & \si{\milli\liter \per \mmHg}   &
		$\EaMaxRV$   & 0.55       & \si{\mmHg \per \milli\liter}   \\
		$\CvnSYS$    & 60.0         & \si{\milli\liter \per \mmHg}   &
		$\VnLA$      & 4.0          & \si{\milli\liter}   \\
		$\CvnPUL$    & 16.0         & \si{\milli\liter \per \mmHg}   &
		$\VnRA$      & 4.0          & \si{\milli\liter}   \\
		$\LarSYS$    & \num{5e-3}         & \si{\mmHg \second\squared \per \milli\liter}  &
		$\VnRV$      & 16.0         & \si{\milli\liter}   \\
		$\LarPUL$    & \num{5e-4}         & \si{\mmHg \second\squared \per \milli\liter}  &
		\multicolumn{3}{l}{\textbf{Cardiac valves}} \\
		$\LvnSYS$    & \num{5e-4}         & \si{\mmHg \second\squared \per \milli\liter}  &
		$\Rmin$      & 0.0075       & \si{\mmHg \second \per \milli\liter}   \\
		$\LvnPUL$    & \num{5e-4}         & \si{\mmHg \second\squared \per \milli\liter}  &
		$\Rmax$      & 75006.2      & \si{\mmHg \second \per \milli\liter}   \\
		\bottomrule
	\end{tabular}
	\caption{Parameters of the circulation model.}
	\label{tab:parameters circulation}
\end{table}

	\end{appendices}

    \clearpage
	\printbibliography

\end{document}